\documentclass[amstex,12pt, amssymb]{article}

\usepackage{mathtext}
\usepackage[cp1251]{inputenc}
\usepackage[T2A]{fontenc}
\usepackage[dvips]{graphicx}
\usepackage{amsmath}
\usepackage{amssymb}
\usepackage{amsxtra}
\usepackage{latexsym}
\usepackage{ifthen}

\newcounter{lemma}[section]

\newcounter{corollary}[section]

\newcounter{remark}[section]

\newcounter{theorem}[section]

\newcounter{proposition}[section]

\newcounter{example}

\numberwithin{equation}{section}

\begin{document}

\markboth{V. DESYATKA, E.~SEVOST'YANOV}{\centerline{ON ISOLATED
SINGULARITIES ...}}

\def\cc{\setcounter{equation}{0}
\setcounter{figure}{0}\setcounter{table}{0}}

\overfullrule=0pt


\author{VICTORIA DESYATKA, EVGENY SEVOST'YANOV}

\title{
{\bf ON ISOLATED SINGULARITIES OF MAPPINGS WITH INVERSE MODULI
INEQUALITIES}}

\date{\today}
\maketitle

\begin{abstract}
We consider open discrete mappings that satisfy the modulus
condition of the inverse Poletsky inequality type. We study the case
when the majorant in it is integrable, or more generally, has finite
averages over infinitesimal spheres. We proved that such mappings
have a continuous extension to an isolated point of the boundary of
some domain without any a priori requirements on the corresponding
mapped domain for integrable majorants. In the case of majorants,
integrable on spheres, we require only the boundedness of the mapped
domain. We do not require any other topological conditions on the
mappings.
\end{abstract}

\bigskip
{\bf 2010 Mathematics Subject Classification: Primary 30C65;
Secondary 31A15, 31B25}

\section{Introduction}

As is known, quasiregular mappings $f:D\rightarrow {\Bbb R}^n,$
$n\geqslant 2,$ satisfy the condition
\begin{equation}\label{eq2}
M(\Gamma)\leqslant K N(f, D) M(f(\Gamma))
\end{equation}
for any family $\Gamma$ of paths $\gamma$ in a domain $D,$ where $M$
is a conformal modulus of families of paths, $K={\rm ess \sup}\,
K_O(x, f),$
$$K_{O}(x,f)\quad =\quad \left\{
\begin{array}{rr}
\frac{\Vert f^{\,\prime}(x)\Vert^n}{|J(x,f)|}, & J(x,f)\ne 0,\\
1, & f^{\,\prime}(x)=0, \\
\infty, & {\rm in\,other\,cases}
\end{array}
\right.\,,$$
$$
\Vert f^{\,\prime}(x)\Vert\,=\,\max\limits_{h\in {\Bbb R}^n
\backslash \{0\}} \frac {|f^{\,\prime}(x)h|}{|h|}\,,\quad
J(x,f)=\det f^{\,\prime}(x)\,,$$
\begin{equation}\label{eq23}
N(y, f, D)\,=\,{\rm card}\,\left\{x\in D: f(x)=y\right\}\,,
\end{equation}
$$N(f, D)\,=\,\sup\limits_{y\in{\Bbb R}^n}\,N(y, f, D)\,,$$
see, e.g., \cite[Theorem~3.2]{MRV$_1$}. The main properties and
behavior of such mappings are currently sufficiently studied. In
particular, quasiregular mappings have a continuous extension to an
isolated point of the boundary of the domain under some conditions
(see, e.g., \cite[Theorem~2.9.III]{Ri}). In this manuscript, we
consider the continuous extension to isolated boundary points of
mappings more general than quasiregular. It is worth noting that a
similar problem was investigated by us in the case of homeomorphisms
and a conformal modulus (see, e.g., \cite[theorem~5.1]{SevSkv$_2$}),
as well as for open discrete mappings $f:D\rightarrow D ^{\,\prime}$
with an additional condition
$$C(x_0, f)\subset
\partial D^{\,\prime}$$
(see, e.g., \cite[Theorem~1]{Sev$_1$}). Here, us usually,
\begin{equation}\label{eq5}
C(x, f):=\{y\in \overline{{\Bbb R}^n}:\exists\,x_k\in D:
x_k\rightarrow x, f(x_k) \rightarrow y, k\rightarrow\infty\}\,.
\end{equation}
Unlike \cite{SevSkv$_2$} and \cite{Sev$_1$}, we consider the
situation of not only the conformal modulus~$M,$ but also the
modulus of arbitrary order $p\geqslant n.$ In addition, we will
consider open discrete mappings {\it without any additional
conditions nor on $C(x_0, f),$ nor on the mapped domain
$D^{\,\prime}.$} In some individual cases, we assume the boundedness
of $D^{\,\prime},$ however, other conditions of the geometry of this
domain are missing.

\medskip
Let us recall some definitions. Everywhere below, the boundary and
the closure of the set $A\subset {\Bbb R}^n$ are denoted by
$\partial A$ and $\overline{A},$ respectively, and they should be
understood in the sense of the extended Euclidean space
$\overline{{\Bbb R}^n}={\Bbb R}^n\cup\{\infty\}.$ Also, everywhere
below $D$ is a domain in~${\Bbb R}^n,$ $n\geqslant 2.$

A Borel function $\rho:{\Bbb R}^n\,\rightarrow [0,\infty] $ is
called {\it admissible} for the family $\Gamma$ of paths $\gamma$ in
${\Bbb R}^n,$ if the relation
\begin{equation}\label{eq1.4}
\int\limits_{\gamma}\rho (x)\, |dx|\geqslant 1
\end{equation}
holds for all (locally rectifiable) paths $ \gamma \in \Gamma.$ In
this case, we write: $\rho \in {\rm adm} \,\Gamma .$ Let $p\geqslant
1,$ then {\it $p$-modulus} of $\Gamma $ is defined by the equality
\begin{equation}\label{eq1.3gl0}
M_p(\Gamma)=\inf\limits_{\rho \in \,{\rm adm}\,\Gamma}
\int\limits_{{\Bbb R}^n} \rho^p (x)\,dm(x)\,.
\end{equation}
We set $M(\Gamma):=M_n(\Gamma).$ Let $y_0\in {\Bbb R}^n,$
$0<r_1<r_2<\infty$ and
\begin{equation}\label{eq1**}
A(y_0, r_1,r_2)=\left\{ y\,\in\,{\Bbb R}^n:
r_1<|y-y_0|<r_2\right\}\,.\end{equation}
Given sets $E,$ $F\subset\overline{{\Bbb R}^n}$ and a domain
$D\subset {\Bbb R}^n$ we denote by $\Gamma(E,F,D)$ the family of all
paths $\gamma:[a,b]\rightarrow \overline{{\Bbb R}^n}$ such that
$\gamma(a)\in E,\gamma(b)\in\,F$ and $\gamma(t)\in D$ for $t \in (a,
b).$ If $f:D\rightarrow {\Bbb R}^n$ is a given mapping, $y_0\in
f(D)$ and $0<r_1<r_2<d_0=\sup\limits_{y\in f(D)}|y-y_0|,$ then by
$\Gamma_f(y_0, r_1, r_2)$ we denote the family of all paths $\gamma$
in $D$ such that $f(\gamma)\in \Gamma(S(y_0, r_1), S(y_0, r_2),
A(y_0,r_1,r_2)).$ Let $Q:{\Bbb R}^n\rightarrow [0, \infty]$ is a
Lebesgue measurable function. We say that {\it $f$ satisfies inverse
Poletsky inequality relative to $p$-modulus} at the point $y_0\in
\overline{f(D)},$ $p\geqslant 1,$ if the ratio
\begin{equation}\label{eq2*A}
M_p(\Gamma_f(y_0, r_1, r_2))\leqslant \int\limits_{f(D)\cap
A(y_0,r_1,r_2)} Q(y)\cdot \eta^p (|y-y_0|)\, dm(y)
\end{equation}
holds for any Lebesgue measurable function $\eta:
(r_1,r_2)\rightarrow [0,\infty ]$ such that
\begin{equation}\label{eqA2}
\int\limits_{r_1}^{r_2}\eta(r)\, dr\geqslant 1\,.
\end{equation}
Note that estimates of the type~(\ref{eq2*A}) are well known and
hold at least for $p=n$ in many classes of mappings (see, e.g.,
\cite[Theorem~3.2]{MRV$_1$}, \cite[Theorem~6.7.II]{Ri} and
\cite[Theorem~8.5]{MRSY}). For $p\ne n$, similar estimates may be
found, e.g., in \cite{GGR}, \cite{MU} and~\cite{U}.

A mapping $f:D\rightarrow {\Bbb R}^n$ is called {\it discrete} if
the image $\{f^{-1}\left(y\right)\}$ of any point $y\,\in\,{\Bbb
R}^n$ consists of isolated points, and {\it open} if the image of
any open set $U\subset D$ is an open set in ${\Bbb R}^n.$

\medskip
Later, in the extended space $\overline{{{\Bbb R}}^n}={{\Bbb
R}}^n\cup\{\infty\}$ we use the {\it spherical (chordal) metric}
$h(x,y)=|\pi(x)-\pi(y)|,$ where $\pi$ is a stereographic projection
of $\overline{{{\Bbb R}}^n}$ onto the sphere
$S^n(\frac{1}{2}e_{n+1},\frac{1}{2})$ in ${{\Bbb R}}^{n+1},$ namely,
$$h(x,\infty)=\frac{1}{\sqrt{1+{|x|}^2}}\,,$$
\begin{equation}\label{eq3C}
\ \ h(x,y)=\frac{|x-y|}{\sqrt{1+{|x|}^2} \sqrt{1+{|y|}^2}}\,, \ \
x\ne \infty\ne y
\end{equation}
(see, e.g., \cite[Definition~12.1]{Va}).
Further, for the sets $A, B\subset \overline{{\Bbb R}^n}$ let us put
\begin{equation}\label{eq1A}
h(A, B)=\inf\limits_{x\in A, y\in B}h(x, y)\,,\quad
h(A)=\sup\limits_{x, y\in A}h(x ,y)\,, \end{equation}
where $h$ is the chordal distance defined in ~(\ref{eq3C}). 
We will put 
%
%
The following statements are true.

\medskip
\begin{theorem}\label{th1}
{\sl\, Let $n\geqslant 2,$ $p\geqslant n,$ let $D$ be a domain in
${\Bbb R}^n,$ $x_0\in D,$ and let $f:D\setminus\{x_0\}\rightarrow
{\Bbb R}^n$ be an open discrete mapping that satisfies the
conditions~(\ref{eq2*A})--(\ref{eqA2}) at any point $y_0\in
\overline{D^{\,\prime}}\setminus\{\infty\},$ where
$D^{\,\prime}:=f(D\setminus\{x_0\}).$

\medskip
If $Q\in L^1(D^{\,\prime}),$ then $f$ has a continuous
extension $\overline{f}:D\rightarrow\overline{{\Bbb R}^n},$ the
continuity of which should be understood in the sense of the chordal
metric $h$ in~(\ref{eq3C}). The extended mapping $\overline{f}$ is
open and discrete in $D.$ Moreover, if $p=n$ and
$\overline{f}(x_0)\ne \infty,$ then there is a neighborhood
$U\subset D$ of the point $x_0$ depending only on $x_0,$ and $C=C(n,
D, x_0)>0$ such that
\begin{equation}\label{eq2C}
|\overline{f}(x)-\overline{f}(x_0)|\leqslant\frac{C_n\cdot (\Vert
Q\Vert_1)^{1/n}}{\log^{1/n}\left(1+\frac{\delta}{2|x-x_0|}\right)}
\end{equation}
for any $x, y\in U,$
where $\Vert Q\Vert_1$ is the norm of the function $Q$ in
$L^1(D^{\,\prime}).$ }
\end{theorem}

\medskip
Under the additional requirement that the mapped domain
$D^{\,\prime}$ is bounded, the conditions on the function $Q$ in
Theorem~\ref{th1} can be somewhat relaxed. The following statement
holds.

\medskip
\begin{theorem}\label{th2}
{\sl\, Let $n\geqslant 2,$ $p\geqslant n,$ let $D$ be a domain in
${\Bbb R}^n,$ $x_0\in D,$ and let $f:D\setminus\{x_0\}\rightarrow
{\Bbb R}^n$ be an open discrete mapping that satisfies
conditions~(\ref{eq2*A})--(\ref{eqA2}) at any point $y_0\in
\partial D^{\,\prime}\setminus\{\infty\},$ where $D^{\,\prime}:=f(D\setminus\{x_0\})$
is a bounded domain.

\medskip
Assume that, in addition, for any $y_0\in
\partial D^{\,\prime}\setminus\{\infty\}$ and any
$0<r_1<r_2<r_0:=\sup\limits_{y\in D^{\,\prime}}|y-y_0|$ there exists
a set $E\subset[r_1, r_2]$ of positive linear Lebesgue measure such
that the function $Q$ is integrable on $S(y_0, r)$ relative to the
$(n-1)$-dimensional Hausdorff measure $\mathcal{H}^{n-1}$ on $S(y_0,
r)$ for any $r\in E.$ Then $f$ has a continuous extension
$\overline{f}:D\rightarrow{\Bbb R}^n,$ the continuity of which
follows to be understood in the sense of the Euclidean metric. The
extension mentioned above is open and discrete in $D.$

In particular, the statement of Theorem~\ref{th2} holds if the
condition on $Q$ is replaced by a simpler one: $Q\in
L^1(D^{\,\prime}).$ }
\end{theorem}

\medskip
It should be noted that the mapping assumptions in Theorem~\ref{th1}
apply to all points $y_0\in
\overline{D^{\,\prime}}\setminus\{\infty\},$ and in
Theorem~\ref{th2} these conditions apply only to boundary points:
$y_0\in \partial D^{\,\prime}\setminus\{\infty\}.$ A similar remark
applies to the function $Q,$ because in Theorem~\ref{th1} we use its
integrability by the Lebesgue measure, and in Theorem~\ref{th2} we
require only its integrability over spheres centered at the boundary
points. In particular, in Theorem~\ref{th2} we do not assume any
behavior of the function $Q$ inside domain $D^{\,\prime}.$

\medskip
\begin{remark}\label{rem2}
Note that, all quasiregular mappings $f:D\setminus\{x_0\}\rightarrow
{\Bbb R}^n$ satisfy the condition
\begin{equation}\label{eq22}
M(\Gamma_f(y_0, r_1, r_2))\leqslant \int\limits_{f(D)\cap
A(y_0,r_1,r_2)} K_O\cdot N(y, f, D\setminus\{x_0\})\cdot \eta^n
(|y-y_0|)\, dm(y) \end{equation}
at any point $y_0\in \overline{f(D)}\setminus\{\infty\}$ with some
constant $K_O=K_O(f)\geqslant 1$ and any Lebesgue measurable
function $\eta: (r_1,r_2)\rightarrow [0,\infty],$ which satisfies
condition~(\ref{eqA2}). Here the multiplicity function $N(y, f,
D\setminus\{x_0\})$ is defined by~(\ref{eq23}).

Indeed, quasiregular mappings satisfy the condition
\begin{equation}\label{eq24}
M(\Gamma_f(y_0, r_1, r_2))\leqslant \int\limits_{f(D)\cap
A(y_0,r_1,r_2)} K_O\cdot N(y, f, D\setminus\{x_0\})\cdot
(\rho^{\,\prime})^n(y)\, dm(y) \end{equation}
for any function $\rho^{\,\prime}\in{\rm adm}\, f(\Gamma_f(y_0, r_1,
r_2)),$ see \cite[Remarks~2.5.II]{Ri}. Now, we put
$\rho^{\,\prime}(y):=\eta(|y-y_0|)$ for $y\in A(y_0,r_1,r_2)\cap
f(D),$ and $\rho^{\,\prime}(y)=0$ otherwise. By the Luzin theorem,
we may assume that the function $\rho^{\,\prime}$ is Borel
measurable (see, e.g., \cite[Section~2.3.6]{Fe}). Then,
by~\cite[Theorem~5.7]{Va} we have that
$$\int\limits_{\gamma_*}\rho^{\,\prime}(y)\,|dy|\geqslant
\int\limits_{r_1}^{r_2}\eta(r)\,dr\geqslant 1$$
for any (locally rectifiable) path $\gamma_*\in \Gamma(S(y_0, r_1),
S(y_0, r_2), A(y_0, r_1, r_2)).$ Substituting the function
$\rho^{\,\prime}$ in~(\ref{eq24}), we obtain the desired
ratio~(\ref{eq22}).
\end{remark}

\medskip
Due to Remark~\ref{rem2} and Theorem~\ref{th2}, we obtain the
following.

\medskip
\begin{corollary}\label{cor1}
{\sl\, Let $n\geqslant 2,$ let $D$ be a domain in ${\Bbb R}^n,$
$x_0\in D,$ and let $f:D\setminus\{x_0\}\rightarrow {\Bbb R}^n$ be a
quasiregular mapping. If $N(y, f, D\setminus\{x_0\})\in
L^1(D^{\,\prime}),$ then $f$ has a continuous extension
$\overline{f}:D\setminus\overline{{\Bbb R}^n},$ the continuity of
which must be understood in the sense of the chordal metric $h$
in~(\ref{eq3C}).}
\end{corollary}

\medskip
The statement of Corollary~\ref{cor1} is not new. It follows from
the Sokhotsky-Weierstrass theorem, which states that, if $x_0$ is an
essential singular point of the mapping $f$, then $N(y, f,
D\setminus\{x_0\})=\infty $ for any $y\in\overline{{\Bbb
R}^n}\setminus E,$ where $E$ is some set of capacity zero (see
\cite[Corollary~2.11.III]{Ri}). So, the condition $N(y, f,
D\setminus\{x_0\})\in L^1(D^{\,\prime})$ which present in
Corollary~\ref{cor1}, a priori cannot to be fulfilled for essential
singular points $x_0$ of $f.$ However, the methodology of proving
Theorems~\ref{th1} and~\ref{th2} (and, therefore,
Corollary~\ref{cor1}) is fundamentally new in the theory of
mappings.

\section{Auxiliary statements} Before starting the
proof of the main ones results, we establish the following
statements.

\medskip
\begin{proposition}\label{pr1}
{\sl\, Let $D$ be a domain in ${\Bbb R}^n,$ $n\geqslant 2,$ $x_0\in
D,$ and let $f:D\setminus\{x_0\}\rightarrow {\Bbb R}^n,$ $n\geqslant
2,$ be an arbitrary mapping. Then $C(x_0, f)$ is a continuum in
$\overline{{\Bbb R}^n}.$ }
\end{proposition}

\medskip
\begin{proof}
Put $r_0=d(z_0, \partial D).$ Let $r_m>0,$ $m\in {\Bbb N},$ be an
arbitrary sequence such that $r_m\rightarrow 0$ as $m\rightarrow
\infty.$ Let $m_0\in {\Bbb N}$ be such that $r_m<r_0$ for all
$m>m_0.$ Note that
$$C(x_0, f)=
\bigcap\limits_{m=1}^{\infty}\overline{f(B(x_0, r_m)\setminus\{x_0\})}\,.$$
It follows that $C(x_0, f)$ is a closed as a countable intersection
of closed sets. In addition, $C(x_0, f)$ is a continuum as
intersection of the decreasing sequence of continua (see
\cite[Theorem~5.II.47.5]{Ku}).~$\Box$
\end{proof}

\medskip Let us first formulate a simple but
very important topological statement, which is repeatedly used later
(see, e.g., \cite[Theorem~1.I.5.46]{Ku}).

\medskip
\begin{proposition}\label{pr2}
{\sl\, Let $A$ be a set in a topological space $X.$ If $C$ is
connected and $C\cap A\ne \varnothing\ne C\setminus A,$ then $C\cap
\partial A\ne\varnothing.$}
\end{proposition}

\medskip
Let $D\subset {\Bbb R}^n,$ $f:D\rightarrow {\Bbb R}^n$ be a discrete
open mapping, $\beta: [a,\,b)\rightarrow {\Bbb R}^n$ be a path, and
$x\in\,f^{\,-1}(\beta(a)).$ A path $\alpha: [a,\,c)\rightarrow D$ is
called a {\it maximal $f$-lifting} of $\beta$ starting at $x,$ if
$(1)\quad \alpha(a)=x\,;$ $(2)\quad f\circ\alpha=\beta|_{[a,\,c)};$
$(3)$\quad for $c<c^{\prime}\leqslant b,$ there is no a path
$\alpha^{\prime}: [a,\,c^{\prime})\rightarrow D$ such that
$\alpha=\alpha^{\prime}|_{[a,\,c)}$ and $f\circ
\alpha^{\,\prime}=\beta|_{[a,\,c^{\prime})}.$ The following
assertion holds (see~\cite[Lemma~3.12]{MRV$_2$}).

\medskip
If $\beta:[a, b)\rightarrow\overline{{\Bbb R}^n}$ is a path and if
$C\subset\overline{{\Bbb R}^n},$ we say that $\beta\rightarrow C$ as
$t\rightarrow b,$ if the spherical distance $h(\beta(t),
C)\rightarrow 0$ as $t\rightarrow b$ (see
\cite[section~3.11]{MRV$_2$}), where $h(\beta(t), C)$ is defined
in~(\ref{eq1A}).

\medskip
\begin{proposition}\label{pr3}
{\sl Let $f:D\rightarrow {\Bbb R}^n,$ $n\geqslant 2,$ be an open
discrete mapping, let $x_0\in D,$ and let $\beta: [a,\,b)\rightarrow
{\Bbb R}^n$ be a path such that $\beta(a)=f(x_0)$ and such that
either $\lim\limits_{t\rightarrow b}\beta(t)$ exists, or
$\beta(t)\rightarrow \partial f(D)$ as $t\rightarrow b.$ Then
$\beta$ has a maximal $f$-lifting $\alpha: [a,\,c)\rightarrow D$
starting at $x_0.$ If $\alpha(t)\rightarrow x_1\in D$ as
$t\rightarrow c,$ then $c=b$ and $f(x_1)=\lim\limits_{t\rightarrow
b}\beta(t).$ Otherwise $\alpha(t)\rightarrow \partial D$ as
$t\rightarrow c.$}
\end{proposition}

\medskip
Everywhere below $\mu(y, f, G)$ denotes the topological index of the
mapping $f$ at the point $y\,\in\,f(G)\setminus f(\partial G)$ with
respect to a domain $G\subset D,$ $\overline{G}\subset D.$ We say
that $f$ is {\it sense-preserving} ({\it sense-reversing}) if the
topological index $\mu(y, f, G)$ satisfies the condition $\mu(y, f,
G)>0$ ($\mu(y, f, G)<0$) for any domain $G\subset D,$
$\overline{G}\subset D,$ and $y\,\in\,f(G)\setminus f(\partial G)$
(see, e.g., item~4, Ch.~I, p.~17 in~\cite{Ri}). Recall that, for any
open discrete mapping $f:D\rightarrow {\Bbb R}^n$ there exists a
domain $G\subset \overline{D}$ such that $\overline{G}\cap
f^{\,-1}(f(x))=\{x\}$ (see \cite[Lemma~4.9.I]{Ri}). Then the
quantity $i(x,f):=\mu(f(x), f, G)$ does not depend on the choice of
such a domain $G,$ and is called the {\it local topological index of
$f$ at $x.$} The following result is true.

\medskip
\begin{proposition}\label{pr4}
{\sl\,Let $n\geqslant 2,$ $D$ be a domain in ${\Bbb R}^n,$ $x_0\in
D,$ and let $f:D\setminus\{x_0\}\rightarrow { \Bbb R}^n$ be an open
discrete mapping. If $f$ has a continuous extension
$\overline{f}:D\rightarrow \overline{{\Bbb R}^n}$ to a point $x_0,$
then $\overline{f}$ is also an open and discrete mapping. }
\end{proposition}

\medskip
\begin{proof}
Without loss of generality, applying an additional inversion
transformation $\psi(x)=\frac{x}{|x|^2},$ $\infty\mapsto 0,$ we may
assume that $\overline{f}(x_0)\ne\infty.$

It is known that discrete open mappings in ${\Bbb R}^n,$ $n\geqslant
2,$ are either sense-preserving or sense-reversing, see, for
example, item~4, Ch.~I in \cite{Ri}. Now, me may assume that $f$ is
sense-preserving. Let us show that, the extended mapping
$\overline{f}$ is sense-preserving, open and discrete. We denote, as
usual, by $B_f\left(D\setminus\{x_0\}\right)$ the branch set of $f$
in $D\setminus\{ x_0\}.$ Similarly, let $B_{\overline{f}}(D)$ be a
set of all branch points of the extended mapping $\overline{f}.$ If
$\overline{f}$ is a local homeomorphism at some neighborhood of
$x_0,$ there is nothing to prove. Let $x_0\,\in\,
B_{\overline{f}}(D).$ By Chernavsky's theorem, ${\rm
dim}\,B_f(D\setminus\{x_0\})=\,{\rm dim
}\,\overline{f}(B_f(D\setminus\{x_0\}))\,\leqslant\,n-2,$ see, e.g.,
Theorem~4.6.I in \cite{Ri}. (Here ${\rm dim}$ denotes the
topological dimension of the set, see \cite{HW}). Now we obtain that
\begin{equation}\label{eq38*!}
{\rm dim}\,f(B_{\overline{f}}(D))\leqslant n-2\,,
\end{equation}
because
$f(B_{\overline{f}}(D))\,=\,f(B_{\overline{f}}(D\setminus\{x_0\}))\bigcup
\left\{\overline{f}(x_0)\right\},$ in addition, the set
$\left\{\overline{f}(x_0)\right\}$ is closed and topological the
dimension of each of the sets $\overline{f}(B_f(D\setminus\{x_0\}))$
and $\left\{\overline{f}(x_0)\right\}$ does not exceed $n-2,$ see
Corollary~1 Ch.~III item~3 in \cite{HW}. Let $G$ be a domain in $D$
such that $\overline{G}\subset D$ and $y\in\, f(G)\setminus
f(\partial G).$ Then, by virtue of (\ref{eq38*!}), there exists a
point $y_0\,\notin\,f(B_{\overline{f}}(D)),$ belonging to the same
component of the set $\overline{{\Bbb R}^n}\setminus f(\partial G)$
as $y.$

\medskip
Due to the fact that the topological index $\mu(y, f , G)$ is a
constant on each connected component of $\overline{{{\Bbb
R}^n}}\setminus f(\partial G)$ (see property~$D_1$ in
\cite[Proposition~4.4.I]{Ri}), we will have:
$$\mu(y, \overline{f}, G)=
\mu(y_0, f, G)=\sum\limits_{x\,\in\,G\cap \{f^{\,-1}(y_0)\}}
i(x,f)>0\,.$$
Therefore, the mapping $f$ is sense-preserving in $D.$ Finally, for
any $y\,\in f(D),$ due to the discreteness of the mapping $f$ in
$D\setminus\{x_0\}, $ the set $\left\{f^{\,-1}(y)\right\}$ is at
most countable and therefore ${\rm
dim}\,\left\{{\overline{f}}^{\,-1}(y)\right\}=0.$ Thus,
by~\cite{TY}, p.~333, the mapping $\overline{f}$ is open and
discrete, as required.~$\Box$
\end{proof}

\medskip
The following statement contains an elementary connection between
the integrability of a function on spheres and its integrability by
the Lebesgue measure. Note that this statement is only one of the
versions of Fubini's theorem (see, e.g.,
\cite[theorem~8.1.III]{Sa}).

\medskip
\begin{proposition}\label{pr5}
{\sl\, Let $D$ be a domain in ${\Bbb R}^n,$ $n\geqslant 2$ $Q\in
L^1(D^{\,\prime}),$ $Q\equiv 0$ on ${\Bbb R}^n\setminus
D^{\,\prime}.$ Then for every $y_0\in
\overline{D^{\,\prime}}\setminus\{\infty\}$ and each
$0<r_1<r_2<r_0:=\sup\limits_{y\in D^{\,\prime}}| y-y_0|$ the
function of $Q$ is integrable on $S(y_0, r)$ for almost any
$r_1\leqslant r\leqslant r_2$ with respect to the
$(n-1)$-dimensional Hausdorff measure $\mathcal{H}^{n-1}$ on $S(y_0,
r).$}
\end{proposition}

\medskip
\begin{proof}
It follows from the conditions of Proposition~\ref{pr5} that $Q\in
L_{\rm loc}^1({\Bbb R}^n).$ Then by Fubini's theorem (see, e.g.,
\cite[Theorem~8.1.III]{Sa}) we obtain that
$$\int\limits_{r_1<|y-y_0|<r_2}Q(y)\,dm(y)=\int\limits_{r_1}^{r_2}
\int\limits_{S(y_0, r)}Q(y)\,d\mathcal{H}^{n-1}(y)dr<\infty\,.$$
This means that $Q$ is integrable on $S(y_0, r)$ for almost all
$r\in [r_1, r_2],$ as required.~$\Box$
\end{proof}

\medskip
Given a domain $D\subset {\Bbb R}^n,$ $n\geqslant 2,$ and a Lebesgue
measurable function $Q:{\Bbb R}^n\rightarrow [0, \infty],$ we denote
by $\ frak{F}_Q(D)$ the family of all open discrete mappings
$f:D\rightarrow {\Bbb R}^n$ such that the ratio~(\ref{eq2*A}) is
fulfilled for $p=n$ and any point $y_0\in f(D).$ The following
statement holds (see, e.g., \cite[Theorem~1.1]{SSD}).

\medskip
\begin{proposition}\label{pr6}
{\sl\, Let $Q\in L^1({\Bbb R}^n).$ Then there is a constant $C_n>0,$
which depends only on the dimension of the space $n,$ such that, for
any $x_0\in D$ and any $r_0>0$ such that $0<2r_0<{\rm dist}\,(x_0,
\partial D),$ the inequality
\begin{equation}\label{eq2CB}
|f(x)-f(x_0)|\leqslant\frac{C_n\cdot (\Vert
Q\Vert_1)^{1/n}}{\log^{1/n}\left(1+\frac{r_0}{|x-x_0|}\right)}
\end{equation}
holds for any $x\in B(x_0, r_0)$ and $f\in \frak{F}_Q(D),$
where $\Vert Q\Vert_1$ is the norm of the function $Q$ in $L^1({\Bbb
R}^n).$ In particular, the family $\frak{F}_Q(D)$ is equicontinuous
in $D.$}
\end{proposition}

\section{Main Lemmas}

The proof of the main results is step-by-step and is based on the
following logic, in line with which we will prove that

\medskip
1) no two different boundary points of the corresponding mapped
domain can belong to the limit set of the mapping $f$ at the point
$x_0,$

\medskip
2) no two distinct inner points of the corresponding mapped domain
can belong to the limit set of the mapping $f$ at the point $x_0,$

3) consequently, no distinct points can belong to the limit set of
the mapping $f$ at the point $x_0.$ There can be only one such
point, which means exactly one.

\medskip
Let us proceed to the implementation of the first stage. The
following statement holds, see~\cite[Lemma~2.2]{SevSkv$_3$},
cf.~\cite[theorem~10.12]{Va}.

\medskip
\begin{lemma}\label{lem2}
{\sl\, Let $D$ be a domain in ${\Bbb R}^n,$ $n\geqslant 2,$ and
$x_0\in D.$ Then, for any $P>0$ and any neighborhood $U$ of the
point $x_0$ there is a neighborhood $V\subset U$ of the same point,
such that the inequality $M(\Gamma(E, F, D))>P$ holds for any
continua $E, F\subset D$ which intersect $\partial U$ and $\partial
V.$}
\end{lemma}

\medskip
\begin{remark}\label{rem1}
Since the modulus of the family of paths passing through a fixed
point is zero (see section~7.9 in \cite{Va}),  the statement of
Lemma~\ref{lem2} remains true for the case when $x_0$ is an isolated
point of the boundary of the domain $D.$
\end{remark}

\medskip
The following lemma indicates that there cannot be more than one
limit point among the boundary points of the mapping $f$ at the
point $x_0$.

\medskip
\begin{lemma}\label{lem1}
{\sl\, Let $D$ be a domain in ${\Bbb R}^n,$ $n\geqslant 2, $ let
$p\geqslant n,$ $x_0\in D$ and let $f: D\setminus\{x_0\}\rightarrow
{\Bbb R}^n$ be an open discrete mapping that satisfies
conditions~(\ref{eq2*A})--(\ref{eqA2}) at least at one finite point
$y_0\in C(x_0, f)\cap \partial D^{\,\prime},$
$D^{\,\prime}:=f(D\setminus\{x_0\}).$ Assume that, for any
$0<r_1<r_2<r_0:=\sup\limits_{y\in D^{\,\prime}}|y-y_0|$ there is a
set $E\subset[r_1, r_2]$ of positive linear Lebesgue measure such
that the function $Q$ is integrable on $S(y_0, r)$ for any $r\in E$
relative to the $(n-1)$-dimensional Hausdorff measure
$\mathcal{H}^{n-1}$ on $S(y_0, r).$

\medskip
Then the set $C(x_0, f)\cap
\partial D^{\,\prime}$ cannot contain
more than one point.}
\end{lemma}

\medskip
\begin{proof} Without loss of generalization, we may consider that
$D$ is a bounded domain.  We will partially use and improve the
approach, performed in the proof of Theorem~1 in \cite{Sev$_2$}, cf.
Theorem~6 in~\cite{Sev$_1$}. Let us prove this assertion by
contradiction. Let $z_1, z_2\in C(x_0, f)\cap
\partial D^{\,\prime},$ $z_1\ne z_2.$
We may assume that $z_1\ne \infty.$ Then we may also find some
sequences $x_m, x_m^{\,\prime}\in D\setminus \{x_0\},$
$m=1,2,\ldots,$ such that $x_m, x_m^{\,\prime}\rightarrow x_0$ as
$m\rightarrow\infty,$ and $y_m:=f(x_m)\rightarrow z_1$ and
$y^{\,\prime}_m:=f(x_m^{\,\prime})\rightarrow z_2$ as
$m\rightarrow\infty.$ Now, we may find some numbers
$\varepsilon_0>0$ and $\varepsilon_1>0$ such that
\begin{equation}\label{eq6}
\overline{B(z_1, \varepsilon_0)}\cap \overline{B_*(z_2,
\varepsilon_1)}=\varnothing\,,
\end{equation}
where $B_*(z_2, \varepsilon_1)=B(z_2, \varepsilon_1)$ for $z_2\ne
\infty$ and $B_*(z_2, \varepsilon_1)=\{y\in \overline{{\Bbb R}^n}:
h(y, \infty)<\varepsilon_1\}$ fir $z_2=\infty.$ Note that, $B(z_1,
\varepsilon_0)$ is a convex set, and $B_*(z_2, \varepsilon_1)$ is
path connected. In this case, points $y_m$ and $z_1$ may be joined
by the segment $I_m(t)=z_1+t(y_m-z_1),$ $t\in (0, 1),$ which
completely belongs to~$B(z_1, \varepsilon_0).$ Similarly, points
$y_m^{\,\prime}$ and $z_2$ may be joined by a path $J_m(t),$ $t\in
[0, 1],$ which belongs to the ''ball'' $B_*(z_2, \varepsilon_1).$

\medskip
By the construction, $|I|\cap\partial D^{\,\prime}\ne\varnothing\ne
|J|\cap\partial D^{\,\prime}.$ Set
$$t_m:=\sup\limits_{t\in[0, 1]:
I_m(t)\in D^{\,\prime}}t\,,\qquad p_m:=\sup\limits_{t\in[0, 1]:
J_m(t)\in D^{\,\prime}}t\,,$$
$$C^m_1:=I|_{[0, t_m)}\,,\qquad C^m_2:=J|_{[0, p_m)}\,.$$
By Proposition~\ref{pr3}, the paths $C^m_1$ and $C^m_2$ have maximal
$f$-liftings $C^{m\,*}_1:[0, c^m_1)\rightarrow D\setminus\{x_0\}$
and $C^{m\,*}_2:[0, c^m_2)\rightarrow D\setminus\{x_0\}$ at the
origin at the points $x_m$ and $x^{\,\prime}_m,$ respectively. Note
that, the case $C^{m\,*}_1(t)\rightarrow z_0$ as $t\rightarrow
c^m_1-0,$ where $z_0\in D\setminus\{x_0\},$ is impossible, because
in this situation, by Proposition~\ref{pr3} we would have that
$c^m_1=t_m$ and $I_m(t)\rightarrow f(z_0)\in D^{\,\prime},$ which
contradicts the definition of $t_m.$ Now, again by
Proposition~\ref{pr3}
\begin{equation}\label{eq1}
h(C^{m\,*}_1(t), \partial (D\setminus\{x_0\}))\rightarrow 0,\qquad
t\rightarrow c^m_1-0\,.
\end{equation}
Similarly, it may be shown that
\begin{equation}\label{eq1G}
h(C^{m\,*}_2(t), \partial (D\setminus\{x_0\}))\rightarrow 0,\qquad
t\rightarrow c^m_2-0\,.
\end{equation}
By~(\ref{eq1}) and~(\ref{eq1G}), the following four cases are
possible:

\medskip
1) $C^{m\,*}_1(t)\rightarrow
\partial D$ as $t\rightarrow c^m_1-0$ and $C^{m\,*}_2(t)\rightarrow
\partial D$ as $t\rightarrow c^m_2-0$ for every $m\in {\Bbb N};$

\medskip
2) $C^{m\,*}_1(t)\rightarrow
\partial D$ as $t\rightarrow c^m_1-0,$ but there is $m_0\in {\Bbb N}$ such that
$C^{m_0\,*}_2(t)\rightarrow \{x_0\}$ as $t\rightarrow c^{m_0}_2-0;$

\medskip
3) there is $m_0\in {\Bbb N}$ such that $C^{m_0\,*}_1(t)\rightarrow
\{x_0\}$ as $t\rightarrow c^{m_0}_1-0,$ in addition,
$C^{m\,*}_2(t)\rightarrow
\partial D$ as $t\rightarrow c^m_2-0$ for every $m\in {\Bbb N};$

\medskip
4) there are $k_0, m_0\in {\Bbb N}$ such that
$C^{k_0\,*}_1(t)\rightarrow \{x_0\}$ as $t\rightarrow c^k_0-0$ and
$C^{m_0\,*}_2(t)\rightarrow \{x_0\}$ as $t\rightarrow c^{m_0}_2-0.$

\medskip
The cases~2) and~3) may be considered similarly, so it is sufficient
to consider the cases~1), 2) and 4).

\medskip
{\bf Case 1).} Let us assume that, $C^{m\,*}_1(t)\rightarrow
\partial D$ as $t\rightarrow c^m_1-0$ and $C^{m\,*}_2(t)\rightarrow
\partial D$ as $t\rightarrow c^m_2-0$ for any $m\in {\Bbb N}.$  Then there is $\delta_0>0$
such that $d(|C^{m\,*}_1|)\geqslant \delta_0>0$ і
$d(|C^{m\,*}_2|)\geqslant \delta_0>0$ for any $m\in {\Bbb N}$ (see
Figure~\ref{fig1}).
\begin{figure}[h]
\centerline{\includegraphics[scale=0.5]{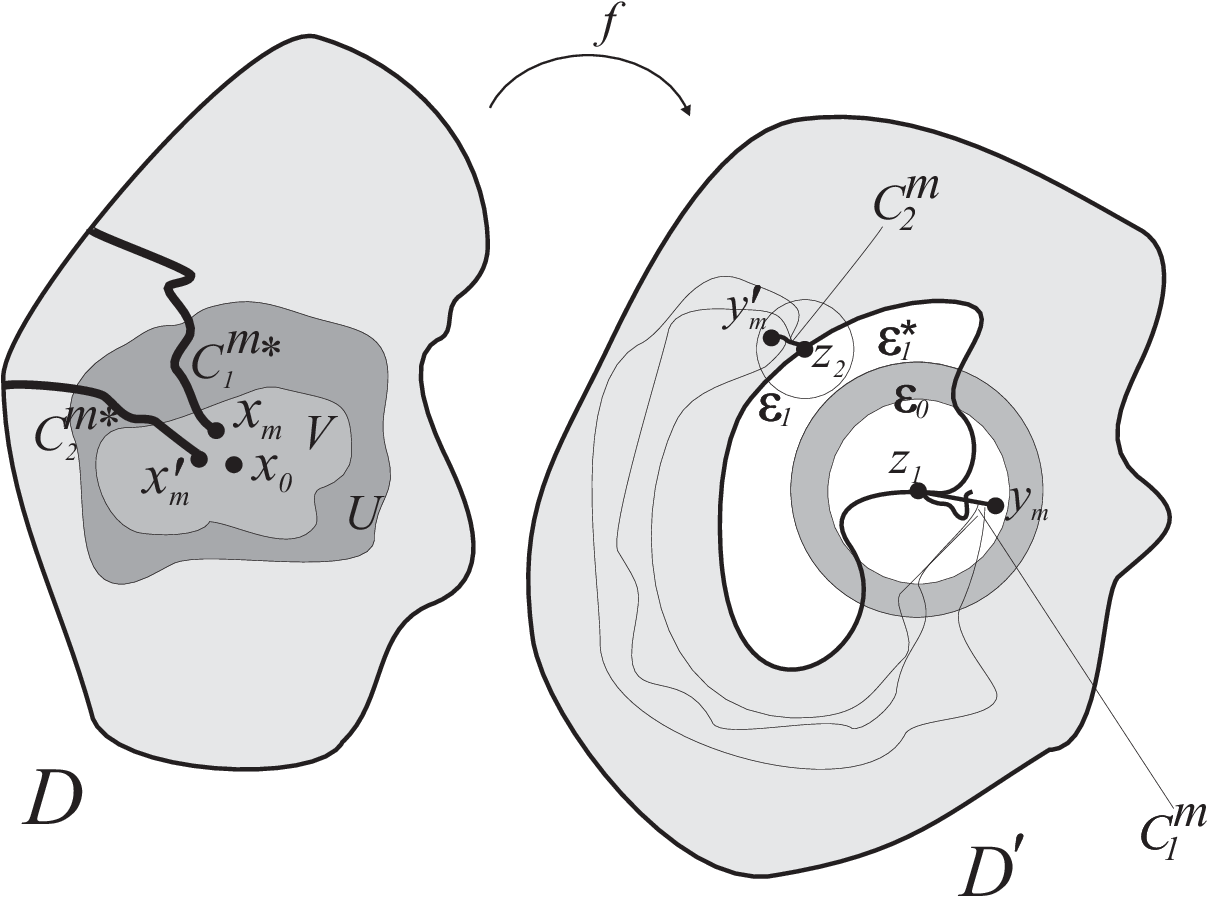}} \caption{To
the proof of Lemma~\ref{lem1}, the case~1)}\label{fig1}
\end{figure}
Put $P>0.$ Let $U:=B(x_0, \delta_0/2),$ and let $V$ be a
neighborhood of $x_0$ which corresponds to Lemma~\ref{lem2} and
Remark~\ref{rem1}. Since by the assumption $x_m, x_m^{\,\prime}\in
D\setminus \{x_0\},$ $m=1,2,\ldots,$ we may find a number $m_0\in
{\Bbb N}$ such that $x_m, x^{\,\prime}_{m}\in V$ for any $m\geqslant
m_0.$
Observe that, for $m\geqslant m_0$
\begin{equation}\label{eq10A}
|C^{m\,*}_1|\cap \partial V\ne\varnothing, \quad |C^{m\,*}_2|\cap
\partial V\ne\varnothing\,.
\end{equation}
Indeed, $x_m\in |C^{m\,*}_1|,$ $x^{\,\prime}_{m}\in |C^{m\,*}_2|,$
and therefore $|C^{m\,*}_1|\cap V\ne\varnothing\ne |C^{m\,*}_2|\cap
V$ for $m\geqslant m_0.$ Besides that, ${\rm diam}\,V\leqslant {\rm
diam}\,U=r_0$ and, since $d(|C^{m\,*}_1|)\geqslant \delta_0>0$ and
$d(|C^{m\,*}_2|)\geqslant \delta_0>0$ form any $m\in {\Bbb N}.$ Now,
by Proposition~\ref{pr2} we obtain the relations~(\ref{eq10A}).
Similarly, we may to prove that
\begin{equation}\label{eq11B}
|C^{m\,*}_1|\cap \partial U\ne\varnothing, \quad |C^{m\,*}_2|\cap
\partial U\ne\varnothing\,.
\end{equation}
Then, by Lemma~\ref{lem2} and Remark~\ref{rem1}
\begin{equation}\label{eq11C}
M(\Gamma(|C^{m\,*}_1|, |C^{m\,*}_2|, D\setminus \{x_0\}))>P\,,
m\geqslant m_0.
\end{equation}
If $p>n,$ since $D$ is bounded, for any $\rho\in {\rm
adm}\,\Gamma(|C^{m\,*}_1|, |C^{m\,*}_2|, D\setminus\{x_0\})$ by the
H\"{o}lder inequality
\begin{equation}\label{eq1B}
M(\Gamma(|C^{m\,*}_1|, |C^{m\,*}_2|, D\setminus\{x_0\}))\leqslant
\int\limits_{D}\rho^n(x)\,dm(x)\leqslant
\left(\int\limits_{D}\rho^p(x)\,dm(x)\right)^{\frac{n}{p}}\cdot
(m(D))^{\frac{p-n}{p}}\,.
\end{equation}
Passing in (\ref{eq1B}) to $\inf$ over all $\rho\in {\rm
adm}\,\Gamma,$ we obtain that
$$M(\Gamma(|C^{m\,*}_1|, |C^{m\,*}_2|,
D\setminus\{x_0\}))\leqslant$$
\begin{equation}\label{eq1E} \leqslant
\int\limits_{D}\rho^n(x)\,dm(x)\leqslant
\left(M_p(\Gamma(|C^{m\,*}_1|, |C^{m\,*}_2|,
D\setminus\{x_0\}))\right)^{\frac{n}{p}}\cdot
(m(D))^{\frac{p-n}{n}}\,.
\end{equation}
It follows from~(\ref{eq11C}) and~(\ref{eq1E}) that
\begin{equation}\label{eq1F}
M_p(\Gamma(|C^{m\,*}_1|, |C^{m\,*}_2|, D\setminus\{x_0\}))\geqslant
P^{\frac{p}{n}}\cdot (m(D))^{\frac{n-p}{n}}\,.
\end{equation}
Let us show that the relation~(\ref{eq11C}) and~(\ref{eq1F}) for
$p=n$ and $p>n,$ respectively, are impossible (in particular, each
of them contradicts with the definition of mapping~$f$
in~(\ref{eq2*A})--(\ref{eqA2})).

\medskip
Since $\overline{B(z_1, \varepsilon_0)}\cap \overline{B_*(z_2,
\varepsilon_1)}=\varnothing,$ we may find
$\varepsilon^*_1>\varepsilon_0$ for which we still have
$\overline{B(z_1, \varepsilon^*_1)}\cap \overline{B_*(z_2,
\varepsilon_1)}=\varnothing.$ Let $\Gamma_*=\Gamma(|C^m_1|, |C^m_2|,
D^{\,\prime}).$ Note that
\begin{equation}\label{eq3D}
\Gamma_*>\Gamma(S(z_1, \varepsilon^*_1), S(z_1, \varepsilon_0),
A(z_1, \varepsilon_0, \varepsilon^*_1))\,.
\end{equation}
Indeed, let $\gamma\in \Gamma_*,$ $\gamma:[a, b]\rightarrow {\Bbb
R}^n.$ Since $\gamma(a)\in |C^m_1|\subset B(z_1, \varepsilon_0)$ and
$\gamma(b)\in |C^m_2|\subset \overline{{\Bbb R}^n}\setminus B(z_1,
\varepsilon_0),$ by Proposition~\ref{pr2} we may find $t_1\in (a,
b)$ such that $\gamma(t_1)\in S(z_1, \varepsilon_0).$ Without loss
of generalization, we may assume that
$|\gamma(t)-z_1|>\varepsilon_0$ for $t>t_1.$ Next, since
$\gamma(t_1)\in S(z_1, \varepsilon^*_1)$ and $\gamma(b)\in
|C^m_2|\subset {\Bbb R}^n\setminus B(z_1, \varepsilon^*_1),$ by
Proposition~\ref{pr2} there is $t_2\in (t_1, b)$ such that
$\gamma(t_2)\in S(z_1, \varepsilon^*_1).$  Without loss of
generalization, we may assume that $|\gamma(t)-z_1|<\varepsilon_0^*$
when $t_1<t<t_2.$ Therefore, $\gamma|_{[t_1, t_2]}$ is a subpath of
$\gamma$ which belongs to~$\Gamma(S(z_1, \varepsilon^*_1), S(z_1,
\varepsilon_0), A(z_1, \varepsilon_0, \varepsilon^*_1)).$ Therefore,
the relation~(\ref{eq3D}) is proved.

\medskip
Let us establish now that
\begin{equation}\label{eq5A}
\Gamma(|C^{m\,*}_1|, |C^{m\,*}_2|, D\setminus\{x_0\})>\Gamma_f(z_1,
\varepsilon_0, \varepsilon^*_1)\,.
\end{equation}
Indeed, if the path $\gamma:[a, b]\rightarrow D\setminus\{x_0\}$
belongs to $\Gamma(|C^{m\,*}_1|, |C^{m\,*}_2|, D\setminus\{x_0\}),$
then $f(\gamma)$ belongs to $D^{\,\prime},$ and $f(\gamma(a))\in
|C^{m\,*}_1|$ and $f(\gamma(b))\in |C^{m\,*}_2|,$ that is,
$f(\gamma)\in \Gamma_*.$ Then according to the above proof and due
to the ratio~(\ref{eq3D}) the path $f(\gamma)$ has a subpath
$f(\gamma)^{\,*}:=f(\gamma)|_{[t_1, t_2]},$ $a\leqslant
t_1<t_2\leqslant b,$ which belongs to the family $\Gamma(S(z_1,
\varepsilon^*_1), S(z_1, \varepsilon_0), A(z_1, \varepsilon_0,
\varepsilon^*_1)).$ Then $\gamma^*:=\gamma|_{[t_1, t_2]}$ is a
subpath of $\gamma$ and it belongs to~$\Gamma_f(z_1, \varepsilon_0,
\varepsilon^*_1),$ as required.

\medskip
In turn, by~(\ref{eq5A}) we have the following:
$$M_p(\Gamma(|C^{m\,*}_1|, |C^{m\,*}_2|, D\setminus\{x_0\}))\leqslant$$
\begin{equation}\label{eq11A}
\leqslant M_p(\Gamma_f(z_1, \varepsilon_0,
\varepsilon^*_1))\leqslant \int\limits_{A} Q(y)\cdot \eta^p
(|y-z_1|)\, dm(y)\,,
\end{equation}
where $A=A(z_1, \varepsilon_0, \varepsilon^*_1)$ is defined
in~(\ref{eq1**}), and $\eta$ is arbitrary Lebesgue measurable
function that satisfies ratio~(\ref{eqA2}) for $r_1:=\varepsilon_0$
and $r_2:=\varepsilon^*_1.$ As above, we use the standard
agreements: $a/\infty=0$ for $a\ne\infty,$ $a/0=\infty$ for $a>0$
and $0\cdot\infty=0$ (see, e.g., \cite[3.I]{Sa}). Put
$\widetilde{Q}(y)=\max\{Q(y), 1\},$
\begin{equation}\label{eq13A}
I=\int\limits_{\varepsilon_0}^{\varepsilon^*_1}\frac{dt}{t^{\frac{n-1}{p-1}}
\widetilde{q}_{z_1}^{1/(p-1)}(t)}\,,
\end{equation}
where $\omega_{n-1}$ denotes the area of the unit sphere ${\Bbb
S}^{n-1}$ in ${\Bbb R}^n,$ and $\widetilde{q}_{z_1}(t)$ is defined
by the relation
\begin{equation}\label{eq12}
\widetilde{q}_{z_1}(r)=\frac{1}{\omega_{n-1}r^{n-1}}\int\limits_{S(z_1,
r)}\widetilde{Q}(y)\,d\mathcal{H}^{n-1}(y)\,. \end{equation}
By the assumption of Lemma~\ref{lem1}, there exists a set $E\subset
[\varepsilon_0, \varepsilon^*_1]$ of positive linear measure such
that $\widetilde{q}_{y_1}(t)$ is finite for all $t\in E.$ Therefore,
$I\ne 0$ in~(\ref{eq13A}). In this case, the function
$\eta_0(t)=\frac{1}{It^{\frac{n-1}{p-1}}
\widetilde{q}_{z_1}^{1/(p-1)}(t)}$ satisfies the
relation~(\ref{eqA2}) for $r_1:=\varepsilon_0$ and
$r_2:=\varepsilon^*_1.$ Substituting this function in the right-hand
side of the inequalities~(\ref{eq11A}) and applying Fubini theorem,
we obtain that
\begin{equation}\label{eq14A}
M_p(\Gamma(|C^{m\,*}_1|, |C^{m\,*}_2|, D\setminus\{x_0\}))\leqslant
\frac{\omega_{n-1}}{I^{p-1}}<\infty\,.
\end{equation}
The relation~(\ref{eq14A}) contradicts with~(\ref{eq11C}) for $p=n$
and~(\ref{eq1F}) for $p>n$. The above contradiction completes the
consideration of the case~1).

\medskip
{\bf Case 2).} $C^{m\,*}_1(t)\rightarrow
\partial D$ as $t\rightarrow c^m_1-0,$ but there is $m_0\in {\Bbb N}$ such that
$C^{m_0\,*}_2(t)\rightarrow \{x_0\}$ as $t\rightarrow c^{m_0}_2-0.$

\medskip
We will show that consideration~2) to a certain extent is reduced to
case~1). Indeed, since $C^{m_0\,*}_2(t)\rightarrow \{x_0\}$ as
$t\rightarrow c^{m_0}_2-0,$ we may find a sequence $t_k\rightarrow
c^{m_0}_2-0$ such that $C^{m_0\,*}_2(t_k)\rightarrow x_0$ as
$k\rightarrow\infty.$ We put $u_k:=C^{m_0\,*}_2(t_k)$ and
$v_k:=f(C^{m_0\,*}_2(t_k)).$ Also let
\begin{equation}\label{eq15}
D^{\,k*}:={C^{m_0\,*}_2}|_{[0, t_k)}\,, \quad
D^k:={C^{m_0\,}_2}|_{[0, t_k)}\,,\quad k=1,2,\ldots\,,
\end{equation}
see Figure~\ref{fig2}.
\begin{figure}[h]
\centerline{\includegraphics[scale=0.5]{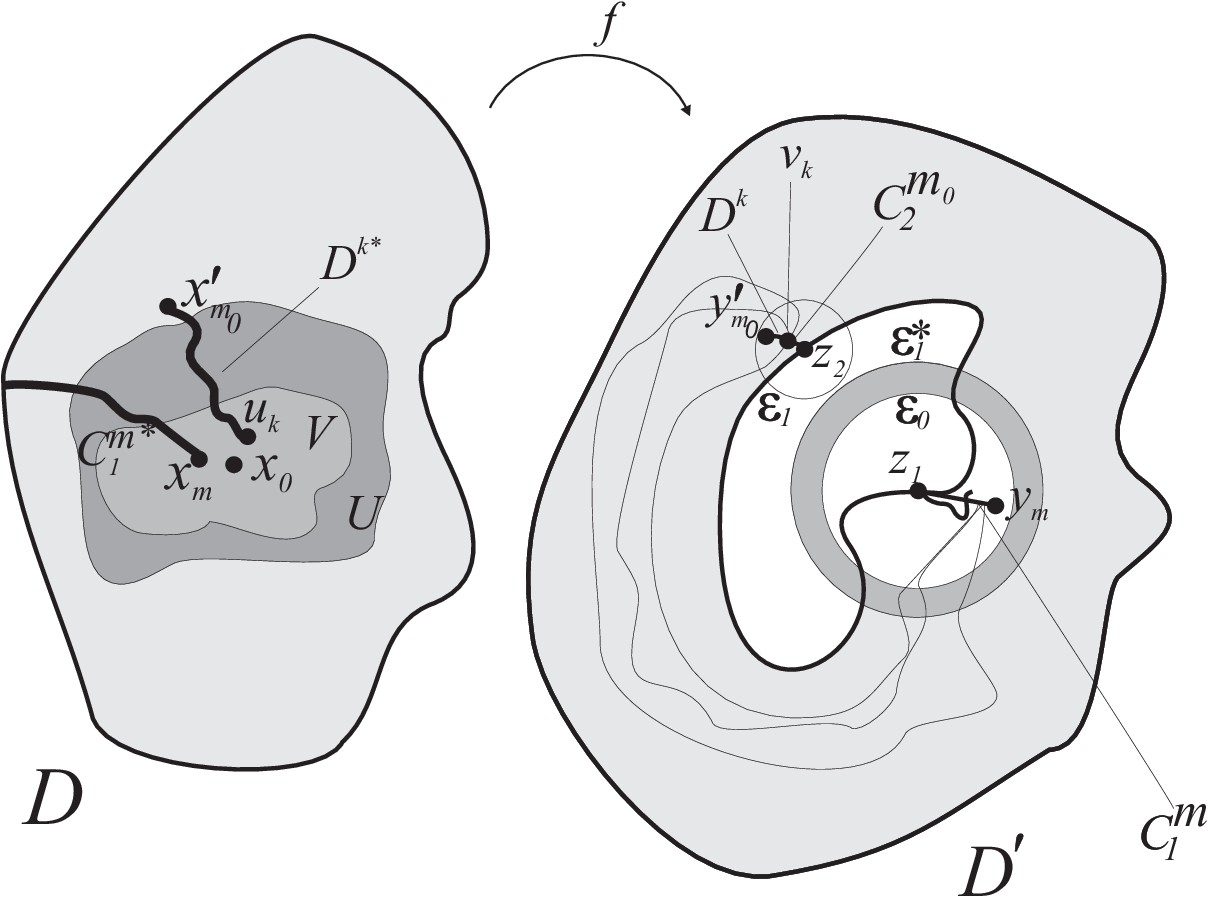}} \caption{To
proof of Lemma~\ref{lem1}, case~2)}\label{fig2}
\end{figure}
Now we reason in the same way as in case~1), where $C^{m\,*}_1$ and
$C^{m}_1$ are constructed in the same way, and the role of
$C^{m\,*}_2$ and $C^{m}_2$ perform $D^{\,m*}$ and $D^{m},$
respectively. Due to the definition $C^{\,m*}_1$ and $D^{\,m*},$
there exists $\delta_0>0$ such that $d(|C^{\,m*}_1|)\geqslant
\delta_0>0$ and $d(|D^{\,m*}|)\geqslant \delta_0>0$ for any
$m=1,2,\ldots .$

\medskip
Let us fix a number $P>0.$ Let $U:=B(x_0, \delta_0/2),$ and let $V$
be a neighborhood of the same point $x_0$ which corresponds to
Lemma~\ref{lem2} and Remark~\ref{rem1}. Reasoning in the same way as
in case~1), we will have
\begin{equation}\label{eq11D}
M(\Gamma(|C^{m\,*}_1|, |D^{m\,*}|, D\setminus \{x_0\}))>P\,,\qquad
m\geqslant m_0
\end{equation}
for $p=n.$ If $p>n,$ since the domain $D$ is bounded,
\begin{equation}\label{eq1H}
M_p(\Gamma(|C^{m\,*}_1|, |D^{m\,*}|, D\setminus\{x_0\}))\geqslant
P^{\frac{p}{n}}\cdot (m(D))^{\frac{n-p}{n}}\,.
\end{equation}
Let us show that the relation~(\ref{eq11D}) and~(\ref{eq1H}) for
$p=n$ and $p>n,$ respectively, are impossible (in particular, each
of them contradicts definition of mapping~$f$
in~(\ref{eq2*A})--(\ref{eqA2})). Really, reasoning in the same way
as in case~1), given that the paths $C^{m}_1$ and $D^{m\,*}$ lie in
disjoint closed balls $\overline{B(z_1, \varepsilon_0)}$ and
$\overline{B_*(z_2, \varepsilon_1)},$ we will have that
\begin{equation}\label{eq14B}
M_p(\Gamma(|C^{m\,*}_1|, |D^{m\,*}|, D\setminus\{x_0\}))\leqslant
\frac{\omega_{n-1}}{I^{p-1}}<\infty\,,
\end{equation}
where $I$ is defined in~(\ref{eq13A}). The relation~(\ref{eq14B})
contradicts~(\ref{eq11D}) for $p=n$ and~(\ref{eq1H}) for $p>n$. The
resulting contradiction completes the consideration case~2).

\medskip
{\bf Case~3)} is considered completely similarly to case~2) (differs
from it by the designations of paths and corresponding points).

\medskip
Finally, consider {\bf case~4)}, which also to some extent reduces
to case~1).

\medskip
Assume that there are $k_0, m_0\in {\Bbb N}$ such that
$C^{k_0\,*}_1(t)\rightarrow \{x_0\}$ as $t\rightarrow c^k_0-0$ and
$C^{m_0\,*}_2(t)\rightarrow \{x_0\}$ as $t\rightarrow c^{m_0}_2-0.$
Let $D^{\,*k}$ and $D^k$ be defined as in~(\ref{eq15}) with
preserving all the relevant notions necessary for them definition.

\medskip
Next, since $C^{k_0\,*}_1(t)\rightarrow \{x_0\}$ as $t\rightarrow
c^{k_0}_2-0,$ there exists a sequence $s_k\rightarrow c^{k_0}_2-0,$
such that $C^{k_0\,*}_1(s_k)\rightarrow \{x_0\}$ as
$k\rightarrow\infty.$ Let us put $w_k:=C^{k_0\,*}_1(s_k)$ and
$\omega_k:=f(C^{k_0\,*}_1(s_k)).$ Let also
\begin{equation}\label{eq15A}
E^{\,k*}:={C^{k_0\,*}_1}|_{[0, s_k)}\,, \quad
E^k:={C^{k_0\,}_1}|_{[0, s_k)}\,,\quad k=1,2,\ldots\,,
\end{equation}
see Figure~\ref{fig3}.
\begin{figure}[h]
\centerline{\includegraphics[scale=0.5]{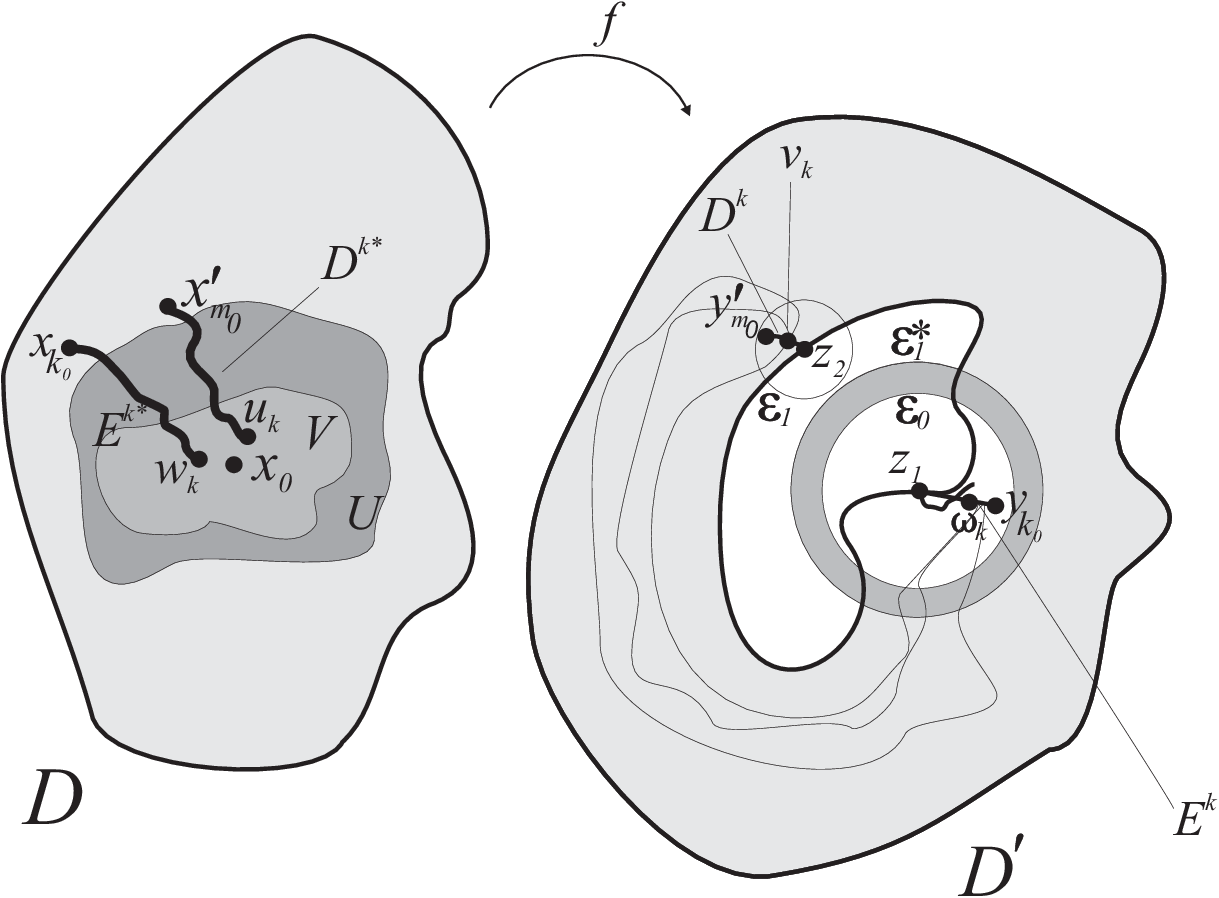}} \caption{To
the proof of Lemma~\ref{lem1}, case~4)}\label{fig3}
\end{figure}
Now we reason in the same way as in case~1), where the role of
$C^{m\,*}_1$ and $C^{m}_1$ perform $E^{\,m*}$ and $E^{m},$
respectively, and the role $C^{\,m*}_2$ and $C^{m}_2$ perform
$D^{\,m*}$ and $D^{m},$ respectively. Due to the definition of
$E^{\,k*}$ and $D^{\,k*},$ there exists $\delta_0>0$ such that
$d(|E^{\,k*}|)\geqslant \delta_0>0$ and $d(|D^{\,k*}|)\geqslant
\delta_0>0$ for all $k\in {\Bbb N}.$

\medskip
Put $P>0.$ Let $U:=B(x_0, \delta_0/2),$ and let $V$ be a
neighborhood of the same point $x_0$ which corresponds to
Lemma~\ref{lem2} and Remark~\ref{rem1}. Reasoning in the same way as
in case~1), we will have
\begin{equation}\label{eq11E}
M(\Gamma(|E^{k\,*}|, |D^{k\,*}|, D\setminus \{x_0\}))>P\,,\quad
k\geqslant K_0
\end{equation}
in the case $p=n.$ If $p>n,$ since $D$ is bounded,
\begin{equation}\label{eq1I}
M_p(\Gamma(|E^{k\,*}|, |D^{k\,*}|, D\setminus\{x_0\}))\geqslant
P^{\frac{p}{n}}\cdot (m(D))^{\frac{n-p}{n}}\,.
\end{equation}
We will show that the relations~(\ref{eq11E}) and~(\ref{eq1I}) for
$p=n$ and $p>n,$, respectively, are impossible (in particular, each
of them contradicts the definition of the mapping~$f$
in~(\ref{eq2*A})--(\ref{eqA2})). Indeed, reasoning in the same way
as in case~1), given that the paths $E^{\,k}$ and $D^{\,k}$ lie in
disjoint balls $\overline{B(z_1, \varepsilon_0)}$ and
$\overline{B_*(z_2, \varepsilon_1)},$ we will have that
\begin{equation}\label{eq14C}
M_p(\Gamma(\Gamma(|E^{k\,*}|, |D^{k\,*}|,
D\setminus\{x_0\}))\leqslant \frac{\omega_{n-1}}{I^{p-1}}<\infty\,,
\end{equation}
where $I$ is defined in~(\ref{eq13A}). The relation~(\ref{eq14C})
contradicts~(\ref{eq11E}) for $p=n$ and~(\ref{eq1I}) for $p>n$. The
resulting contradiction completes the consideration case~4). The
lemma is completely proved.~$\Box$
\end{proof}

\medskip
From Lemma~\ref{lem1} we immediately have the following consequence
(previously it was obtained by us as the main result of the
corresponding paper~\cite[Theorem~1]{Sev$_2$}; see
also~\cite[Theorem~6]{Sev$_1$}).

\medskip
\begin{corollary}\label{cor2}
{\sl\, Let $n\geqslant 2,$ $p\geqslant n,$ $let D$ be a domain in
${\Bbb R}^n,$ $x_0\in D,$ and let $f:D\setminus\{x_0\}
D^{\,\prime},$ $f(D\setminus\{x_0\})=D^{\,\prime}\subset {\Bbb
R}^n,$ be an open discrete mapping that satisfies
conditions~(\ref{eq2*A})--(\ref{eqA2}) at least in one finite point
$y_0\in C(x_0, f).$ Assume that, $C(x_0, f)\subset D^{\,\prime},$
and that for $y_0$ mentioned above and any
$0<r_1<r_2<r_0:=\sup\limits_{y\in D^{\,\prime}}|y-y_0|$ there is a
set $E\subset[r_1, r_2]$ of positive linear Lebesgue measure such
that the function $Q$ is integrable on $S(y_0, r)$ for each $r\in E$
relative to the $(n-1)$-dimensional Hausdorff measure
$\mathcal{H}^{n-1}$ on $S(y_0, r).$

Then $f$ has a continuous extension
$\overline{f}:D\setminus\{x_0\}\rightarrow\overline{{\Bbb R}^n},$
the continuity of which should be understood in the sense of the
chordal metric $h.$ The specified extension is open and discrete in
$D.$

In particular, the statement of Corollary~\ref{cor2} is valid if the
condition on $Q$ is replaced by a simpler one: $Q\in
L^1(D^{\,\prime}).$ }
\end{corollary}

\medskip
\begin{proof}
Due to the compactness of the space $\overline{{\Bbb R}^n},$ the set
$C(x_0, f)$ is not empty. Therefore, the existence of a continuous
extension of the mapping $f$ directly follows from Lemma~\ref{lem1}
taking into account Proposition~\ref{pr5}. The openness and
discreteness of the extended mapping $\overline{f}$ in $D$ follows
from Proposition~\ref{pr4}.~$\Box$
\end{proof}

\medskip
Thus, we have shown that boundary points cannot be cluster points
for the mapping $f$ if there are two or more such points. The
following two lemmas, which are very similar in their formulations
and proofs, contain the assertion that the inner points of the
mapped domain also have the same property: the situation when there
are at least two of them is excluded for mappings that satisfy
conditions~(\ref{eq2*A})--(\ref{eqA2}) under appropriate conditions
on~$Q.$

\medskip
\begin{lemma}\label{lem3}
{\sl\, Let $D$ be a domain in ${\Bbb R}^n,$ $n\geqslant 2,$ $x_0\in
D,$ let $p\geqslant n$ and let $f:D\setminus\{x_0\}\rightarrow
D^{\,\prime},$ $D^{\,\prime}:=f(D\setminus\{x_0\}),$ be an open
discrete mapping that satisfies
conditions~(\ref{eq2*A})--(\ref{eqA2}) for any $y_0\in
\partial D^{\,\prime}.$ Assume that, the domain $D^{\,\prime}$ is bounded and that, for
of any $y_0\in\partial D^{\,\prime}$ and each
$0<r_1<r_2<r_0:=\sup\limits_{y\in D^{\,\prime}}|y-y_0|$ there is a
Lebesgue measurable set $E\subset[r_1, r_2]$ of positive linear
measure such that the function $Q$ is integrable on $S(y_0, r)$ with
respect to $(n-1)$-dimensional Hausdorff measure $\mathcal{H}^{n-1}$
on $S(y_0, r)$ for all $r\in E.$

\medskip
Then the set $C(x_0, f)\cap D^{\,\prime}$ cannot contain more one
point.}
\end{lemma}

\medskip
\begin{proof}
We modify the approach first developed by us
in~\cite[Theorem~1.5]{SevSkv$_2$}, see also~\cite{SevSkv$_1$},
\cite{SevSkv$_3$}, \cite{SevSkv$_4$} and \cite{SSD}.

Let us prove the lemma from the opposite. Suppose the statement
mentioned above is not true. Then, there are at least two sequences
$x_m, y_m\in D\setminus\{x_0\},$ $x_m, y_m\rightarrow x_0$ as
$m\rightarrow\infty$ such that $f(x_m)\rightarrow a$ and
$f(y_m)\rightarrow b$ as $m\rightarrow\infty,$ where $a, b\in
D^{\,\prime}$ and $a\ne b.$ We may assume that $f(x_m)\in B(a,
\varepsilon_0)\subset D$ and $f(y_m)\in B(b, \varepsilon_0)\subset
D$ for all $m\in {\Bbb N},$ where $\varepsilon_0>0$ is such that
\begin{equation}\label{eq13***}
\overline{B(a, \varepsilon_0)}\cap \overline{B(b,
\varepsilon_0)}=\varnothing\,.
\end{equation}
In particular, it follows from ~(\ref{eq13***}) that
\begin{equation}\label{eq13B}
\varepsilon_0<\frac{|a-b|}{2}\,.
\end{equation}
Consider the straight line
$$r=r(t)=a+(b-a)t\,,\quad-\infty<t<\infty\,,$$
which passes through points $a$ and $b,$ see Figure~\ref{fig4}.
\underline{}\begin{figure}[h]
\centerline{\includegraphics[scale=0.5]{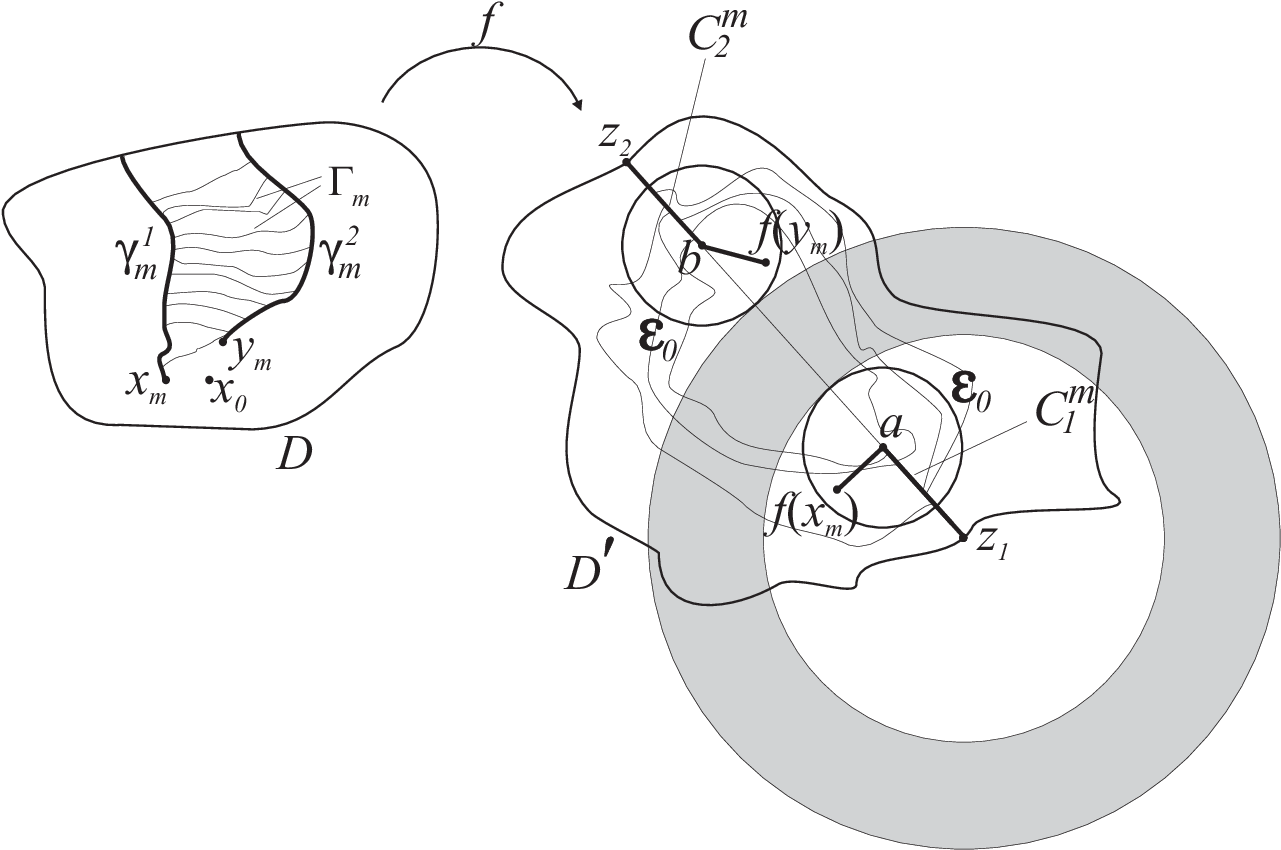}} \caption{To
the proof of Lemma~\ref{lem3}, case~1) }\label{fig4}
\end{figure}
Since $D^{\,\prime}$ is bounded, there are points $z_1$ and $z_2$ of
the intersection of the line $r=r(t)$ with the boundary of the
domain $D^{\,\prime},$ see Proposition~\ref{pr2}. We may assume
that, the points $z_1$ and $z_2$ located on different sides of the
segment $\overline{ab}.$ Let us join the point $f(x_m)$ by the
segment with point $a,$ and point $f(y_m)$ by the segment with point
$b.$ We denote by $C^m_1$ the broken line $\overline{f(x_m)az_1}$
without the point $z_1,$ and, similarly, by $C^m_2$ the broken line
$\overline{f(y_m)bz_2}$ without the point $z_2$ on it. We may assume
that, these broken lines belong completely to $D^{\,\prime},$ i.e.,
all points of segments $\overline{az_1}$ and $\overline{bz_2}$ are
inner points $D^{\,\prime}$ excluding their ends. Let us
parameterize these broken lines as follows:
$$C^m_1:[0, 1)\rightarrow D^{\,\prime}, \quad C^m_1(0)=f(x_m), \quad C^m_1(1)=f(z_1)\,,$$
$$C^m_2:[0, 1)\rightarrow D^{\,\prime}, \quad C^m_2(0)=f(y_m), \quad C^m_2(1)=f(z_2)\,.$$

\medskip
Let us establish that
\begin{equation}\label{eq16}
|C_1^m|\subset B(z_1, |a-z_1|+\varepsilon_0)\,.
\end{equation}
Indeed, if $x\in |C_1^m|,$ then either $x\in \overline{az_1}$ (where
$\overline{az_1}$ is the segment joining the points $a$ and $z_1$),
or $x\in \overline{f(x_m)a}\subset B(a, \varepsilon_0).$ If $x\in
\overline{az_1},$ then, obviously, $|x-z_1|\leqslant
|a-z_1|<|a-z_1|+\varepsilon_0;$ therefore, $x\in B(z_1,
|a-z_1|+\varepsilon_0|).$ Now, if $x\in \overline{f(x_m)a}\subset
B(a, \varepsilon_0),$ then
$|x-a|<\varepsilon_0<|a-z_1|+\varepsilon_0|.$ Therefore also $x\in
B(z_1, |a-z_1|+\varepsilon_0|).$ Inclusion~(\ref{eq16}) is
established.

\medskip
Let us establish now that
\begin{equation}\label{eq17}
|C_2^m|\subset D^{\,\prime}\setminus \overline{B(z_1,
|b-z_1|-\varepsilon_0)}\,.
\end{equation}
Let $x\in |C_2^m|.$ Then either $x\in \overline{bz_2}$ (where
$\overline{bz_2}$ is the segment joining the points $b$ and $z_2$),
or $x\in \overline{f(y_m)b}\subset B(b, \varepsilon_0).$

Assume that $x\in \overline{bz_2}.$ By the construction,
$x=a+(b-a)t_1,$ $t_1\geqslant 1,$ and $z_1=a+(b-a)t_2,$
$t_2\leqslant 0.$ Then
$$|x-z_1|=|a+(b-a)t_1-a-(b-a)t_2|=|(b-a)(t_1-t_2)|\geqslant$$
$$\geqslant |b-a|(1-t_2)= |z_1-b|\geqslant
|b-z_1|-\varepsilon_0\,,$$
therefore $x\in D^{\,\prime}\setminus \overline{B(z_1,
|b-z_1|-\varepsilon_0|)}.$

Now, let $x\in \overline{f(y_m)b}\subset B(b, \varepsilon_0).$ Then
by the triangle inequality
$$|x-z_1|\geqslant |z_1-b|-|b-x|\geqslant |b-z_1|-\varepsilon_0\,.$$
The relation~(\ref{eq17}) is also established.

\medskip
Let $\gamma_m^1:[0, c^1_m)\rightarrow D\setminus\{x_0\},$
$0<c^1_m\leqslant 1,$ be a maximal $f$-lifting of the path $C_1^{m}$
starting at the point $x_m,$ which exists by Proposition~\ref{pr3}.
Let us to prove that the situation $\gamma_m^1\rightarrow
\omega_1\in D\setminus\{x_0\}$ as $t\rightarrow c^1_m$ is not
possible. Indeed, in the contrary case, by Proposition~\ref{pr3} we
obtain that $c^1_m=1$ and $f(\omega_1)=\lim\limits_{t\rightarrow
1-0}C_1^m(t).$ Then, from one side, due to the openness of the
mapping $f,$ $f(\omega_1)\in D^{\,\prime},$ and on the other hand,
$f(\omega_1)=z_1$ by definition of $C_1^m(t).$ Since $z_1\not\in
D^{\,\prime},$ we obtain a contradiction. So, the case
$\gamma_m^1(t)\rightarrow \omega_1\in D\setminus\{x_0\}$ as
$t\rightarrow c^1_m$ is impossible. Again by Proposition~\ref{pr3}
\begin{equation}\label{eq3}
h(\gamma_m^1(t), \partial (D\setminus\{x_0\}))\rightarrow 0
\end{equation}
as $t\rightarrow c^1_m-0.$ Similarly, let $\gamma_m^2:[0,
c^2_m)\rightarrow D\setminus\{x_0\},$ $0<c^2_m\leqslant 1,$ be a
maximal $f$-lifting of the curve $C_2^m$ starting at the point
$y_m,$ which exists by Proposition~\ref{pr3}. Similarly
to~(\ref{eq3}), we obtain that
\begin{equation}\label{eq4}
h(\gamma_m^2(t), \partial D\setminus\{x_0\})\rightarrow 0
\end{equation}
as $t\rightarrow c^2_m-0.$

\medskip
Due to~(\ref{eq3}) and~(\ref{eq4}), the following four cases are
possible:

\medskip
1) for any $m\in {\Bbb N},$ $h(\gamma_{m}^1(t), \partial
D)\rightarrow 0$ as $t\rightarrow c^1_{m}-0$ and $h(\gamma_{m}^2(t),
\partial D) \rightarrow
0$ as $t\rightarrow c^2_{m}-0;$

\medskip
2) for any $m\in {\Bbb N},$ $h(\gamma_{m}^1(t), \partial
D)\rightarrow 0$ as $t\rightarrow c^1_{m}-0,$ but there exists
$m_0\in {\Bbb N}$ such that $h(\gamma_{m_0}^2(t), x_0)\rightarrow 0$
as $t\rightarrow c^2_{m_0}-0;$

\medskip
3) there exists $k_0\in {\Bbb N}$ such that $h(\gamma_{k_0}^1(t),
x_0)\rightarrow 0$ as $t\rightarrow c^1_{k_0}-0,$ in addition, for
all $m\in {\Bbb N},$ $h(\gamma_{m}^2(t),
\partial D)\rightarrow 0$ as $t\rightarrow c^2_{m}-0;$

\medskip
4) there are $k_0, m_0\in {\Bbb N}$ such that $h(\gamma_{k_0}^1(t),
x_0)\rightarrow 0$ as $t\rightarrow c^1_{k_0}-0$ and
$h(\gamma_{m_0}^2(t), x_0)\rightarrow 0$ as $t\rightarrow
c^2_{m_0}-0.$

\medskip
Let us consider each case separately. Note that case~2) is
<<symmetric>> to case~3). Therefore, of these two cases, it is
enough to consider only case~2).

\medskip
Consider {\bf case~1)}: for any $m\in {\Bbb N},$ $h(\gamma_{m}^1(t),
\partial D)\rightarrow 0$ as $t\rightarrow c^1_{m}-0$ and
$h(\gamma_{m}^2(t),
\partial D) \rightarrow
0$ as $t\rightarrow c^2_{m}-0$ (see Figure~\ref{fig4}).

Then there exists $\delta_0>0$ such that $d(|\gamma^1_m|)\geqslant
\delta_0>0$ and $d(|\gamma^2_m|)\geqslant \delta_0>0$ for any $m\in
{\Bbb N}$ (see Figure~\ref{fig4}). Fix $P>0.$ Let $U:=B(x_0,
\delta_0/2),$ and let $V$ be a neighborhood of the same point $x_0,$
which corresponds to Lemma~\ref{lem2} and Remark~\ref{rem1}. Since,
by the assumption, $x_m, y_m\in D\setminus \{x_0\},$ $m=1,2,\ldots,$
there exists a number $m_0\in {\Bbb N}$ such that $x_m, y_m\in V$
when of all $m\geqslant m_0.$
Note that, for $m\geqslant m_0,$
\begin{equation}\label{eq10B}
|\gamma^1_m|\cap \partial V\ne\varnothing, \quad |\gamma^2_m|\cap
\partial V\ne\varnothing\,.
\end{equation}
Indeed, $x_m\in |\gamma^1_m|,$ $y_m\in |\gamma^2_m|,$ so
$|\gamma^1_m|\cap V\ne\varnothing\ne |\gamma^2_m|\cap V$ for
$m\geqslant m_0.$ In addition, ${\rm diam}\,V\leqslant {\rm
diam}\,U=r_0.$ Therefore, since $d(|\gamma^1_m|)\geqslant
\delta_0>0$ and $d(|\gamma^2_m|)\geqslant \delta_0>0$ for $m\in
{\Bbb N},$ by~Proposition~\ref{pr2} we obtain the
relation~(\ref{eq10B}). Similarly, we may prove that
\begin{equation}\label{eq11F}
|\gamma^1_m|\cap \partial U\ne\varnothing, \quad |\gamma^2_m|\cap
\partial U\ne\varnothing\,.
\end{equation}
Then, by Lemma~\ref{lem2} and Remark~\ref{rem1}
\begin{equation}\label{eq11G}
M(\Gamma(|\gamma^1_m|, |\gamma^2_m|, D\setminus \{x_0\}))>P\,,\qquad
m\geqslant m_0.
\end{equation}
If $p>n,$ then, by H\"{o}lder inequality, for every $\rho\in {\rm
adm}\,\Gamma(|\gamma^1_m|, |\gamma^2_m|, D\setminus\{x_0\})$ we
obtain that
\begin{equation}\label{eq1J}
M(\Gamma(|\gamma^1_m|, |\gamma^2_m|, D\setminus\{x_0\}))\leqslant
\int\limits_{D}\rho^n(x)\,dm(x)\leqslant
\left(\int\limits_{D}\rho^p(x)\,dm(x)\right)^{\frac{n}{p}}\cdot
(m(D))^{\frac{p-n}{p}}\,.
\end{equation}
Letting in~(\ref{eq1J}) to $\inf$ over all $\rho\in {\rm
adm}\,\Gamma,$ we have that
$$M(\Gamma(|\gamma^1_m|, |\gamma^2_m|,
D\setminus\{x_0\}))\leqslant$$
\begin{equation}\label{eq1K} \leqslant
\int\limits_{D}\rho^n(x)\,dm(x)\leqslant
\left(M_p(\Gamma(|\gamma^1_m|, |\gamma^2_m|,
D\setminus\{x_0\}))\right)^{\frac{n}{p}}\cdot
(m(D))^{\frac{p-n}{n}}\,.
\end{equation}
By~(\ref{eq11G}) and~(\ref{eq1K}) it follows that
\begin{equation}\label{eq1L}
M_p(\Gamma(|\gamma^1_m|, |\gamma^2_m|, D\setminus\{x_0\}))\geqslant
P^{\frac{p}{n}}\cdot (m(D))^{\frac{n-p}{n}}\,.
\end{equation}
Let us to show that the relations~(\ref{eq11G}) and~(\ref{eq1L}) for
$p=n$ and $p>n,$ respectively, a re impossible (in particular, each
of them contradicts with the definition of the mapping~$f$
in~(\ref{eq2*A})--(\ref{eqA2})).

\medskip
Observe that,
\begin{equation}\label{eq3E}
\Gamma(|C^m_1|, |C^m_2|, D^{\,\prime})>\Gamma(S(z_1, r_1|), S(z_1,
r_2), A(z_1, r_1, r_2)\,,
\end{equation}
where
\begin{equation}\label{eq20}
r_1=|a-z_1|+\varepsilon_0\,, \qquad r_2=|b-z_1|-\varepsilon_0\,.
\end{equation}
First of all, let us to precise that $r_2>r_1.$ Indeed, if
$z_1=a+(b-a)t_2$ for some $t_2\leqslant 0,$ then
$$r_2-r_1=|b-z_1|-\varepsilon_0-|a-z_1|-\varepsilon_0=$$
$$=|b-a-(b-a)t_2|-|a-a-(b-a)t_2|-2\varepsilon_0=$$$$=|b-a|(1-t_2)+|b-a|t_2-2\varepsilon_0=|b-a|-
2\varepsilon_0\geqslant 0\,,$$
see~(\ref{eq13B}).

Now, let $\gamma\in \Gamma(|C^m_1|, |C^m_2|, D^{\,\prime}),$
$\gamma:[0, 1]\rightarrow {\Bbb R}^n.$ Since by~(\ref{eq16})
and~(\ref{eq17}) $\gamma(0)\in |C^m_1|\subset |C_1^m|\subset B(z_1,
r_1)$ and $\gamma(1)\in |C^m_2|\subset D^{\,\prime}\setminus
\overline{B(z_1, r_2)},$ by Proposition~\ref{pr2} there is $t_1\in
(0, 1)$ such that $\gamma(t_1)\in S(z_1, r_1).$ Without loss of
restriction, we may assume that $|\gamma(t)-z_1|>r_1$ при $t>t_1.$
Next, since $\gamma(t_1)\in S(z_1, r_1)$ and $\gamma(1)\in
|C^m_2|\subset {\Bbb R}^n\setminus B(z_1, r_2),$ due to
Proposition~\ref{pr2} there is $t_2\in (t_1, 1)$ such that
$\gamma(t_2)\in S(z_1, r_2).$  Without loss of restriction, we may
assume that $|\gamma(t)-z_1|<r_2$ for $t_1<t<t_2.$ Thus,
$\gamma|_{[t_1, t_2]}$ is a subpath of $\gamma$ which belongs
to~$\Gamma(S(z_1, r_1), S(z_1, r_2), A(z_1, r_1, r_2)).$ Thus, the
relation~(\ref{eq3E}) is proved.

\medskip
Let us to show now, that
\begin{equation}\label{eq5C}
\Gamma(|\gamma^1_m|, |\gamma^2_m|, D\setminus\{x_0\})>\Gamma_f(z_1,
r_1, r_2)\,,
\end{equation}
where, as above, $r_1$ and $r_2$ are defined in~(\ref{eq20}).
Indeed, if $\gamma:[0, 1]\rightarrow D\setminus\{x_0\}$ belongs to
$\Gamma(|\gamma^1_m|, |\gamma^2_m|, D\setminus\{x_0\}),$ then
$f(\gamma)$ belongs to $D^{\,\prime},$ while $f(\gamma(0))\in
|C_1^{m}|$ and $f(\gamma(1))\in |C_2^{m}|.$ Thus, $f(\gamma)\in
\Gamma_*.$ Then, by the proving above and due to the
relation~(\ref{eq3E}) the path $f(\gamma)$ has a subpath
$f(\gamma)^{\,*}:=f(\gamma)|_{[t_1, t_2]},$ $0\leqslant
t_1<t_2\leqslant 1,$ which belongs to the family $\Gamma(S(z_1,
r_2), S(z_1, r_1), A(z_1, r_1, r_2)).$ Then
$\gamma^*:=\gamma|_{[t_1, t_2]}$ is a subpath of $\gamma$ which
belongs to~$\Gamma_f(z_1, r_1, r_2),$ as required.

\medskip
In turn, from~(\ref{eq5C}) it follows that
$$M_p(\Gamma(|\gamma^1_m|, |\gamma^2_m|, D\setminus\{x_0\}))\leqslant$$
\begin{equation}\label{eq11H}
\leqslant M_p(\Gamma_f(z_1, r_1, r_2))\leqslant \int\limits_{A}
Q(y)\cdot \eta^p (|y-z_1|)\, dm(y)\,,
\end{equation}
where $A=A(z_1, r_1, r_2)$ and $\eta$ is any nonnegative Lebesgue
measurable function satisfying the relation~(\ref{eqA2}). As above,
we use the rules $a/\infty=0$ for $a\ne\infty,$ $a/0=\infty$ for
$a>0$ and $0\cdot\infty=0$ (see, e.g., \cite[3.I]{Sa}). Set
$\widetilde{Q}(y)=\max\{Q(y), 1\},$
\begin{equation}\label{eq13C}
I=\int\limits_{r_1}^{r_2}\frac{dt}{t^{\frac{n-1}{p-1}}
\widetilde{q}_{z_1}^{1/(p-1)}(t)}\,,
\end{equation}
where $\omega_{n-1}$ denotes the area of the unit sphere ${\Bbb
S}^{n-1}$ in ${\Bbb R}^n,$ and $\widetilde{q}_{z_1}(t)$ is defined
by the relation
\begin{equation}\label{eq12B}
\widetilde{q}_{z_1}(r)=\frac{1}{\omega_{n-1}r^{n-1}}\int\limits_{S(z_1,
r)}\widetilde{Q}(y)\,d\mathcal{H}^{n-1}(y)\,. \end{equation}
By the assumption, there is a set $E\subset [r_1, r_2]$ of positive
linear Lebesgue measure such that $\widetilde{q}_{z_1}(t)$ is finite
for any $t\in E.$ Thus, $I\ne 0$ in~(\ref{eq13A}). In this case, the
function $\eta_0(t)=\frac{1}{It^{\frac{n-1}{p-1}}
\widetilde{q}_{z_1}^{1/(p-1)}(t)}$ satisfies the
relation~(\ref{eqA2}). Substituting this function in the right-hand
part of~(\ref{eq11H}) and using the Fubini theorem, we obtain that
\begin{equation}\label{eq14D}
M_p(\Gamma(|\gamma^1_m|, |\gamma^2_m|, D\setminus\{x_0\}))\leqslant
\frac{\omega_{n-1}}{I^{p-1}}<\infty\,.
\end{equation}
The relation~(\ref{eq14D}) contradicts with~(\ref{eq11G}) for $p=n$
and~(\ref{eq1L}) for $p>n$. The obtained contradiction completes the
consideration of the case~1).

Consider the {\bf case 2):} for any $m\in {\Bbb N},$
$h(\gamma_{m}^1(t),
\partial D)\rightarrow 0$ as $t\rightarrow c^1_{m}-0,$ however, there exists
$m_0\in {\Bbb N}$ such that $h(\gamma_{m_0}^2(t), x_0)\rightarrow 0$
as $t\rightarrow c^2_{m_0}-0.$

Since $\gamma_{m_0}^2(t)\rightarrow \{x_0\}$ as $t\rightarrow
c^{m_0}_2-0,$ there is a sequence $t_k\rightarrow c^{m_0}_2-0$ such
that $\gamma_{m_0}^2(t_k)\rightarrow x_0.$ Set
$u_k:=\gamma_{m_0}^2(t_k)$ and $v_k:=f(\gamma_{m_0}^2(t_k)).$ Let
also
\begin{equation}\label{eq15B}
D^{\,k*}:={\gamma_{m_0}^2}|_{[0, t_k]}\,, \quad
D^k:={C^{m_0\,}_2}|_{[0, t_k]}\,,\quad k=1,2,\ldots\,,
\end{equation}
see Figure~\ref{fig5}.
\begin{figure}[h]
\centerline{\includegraphics[scale=0.5]{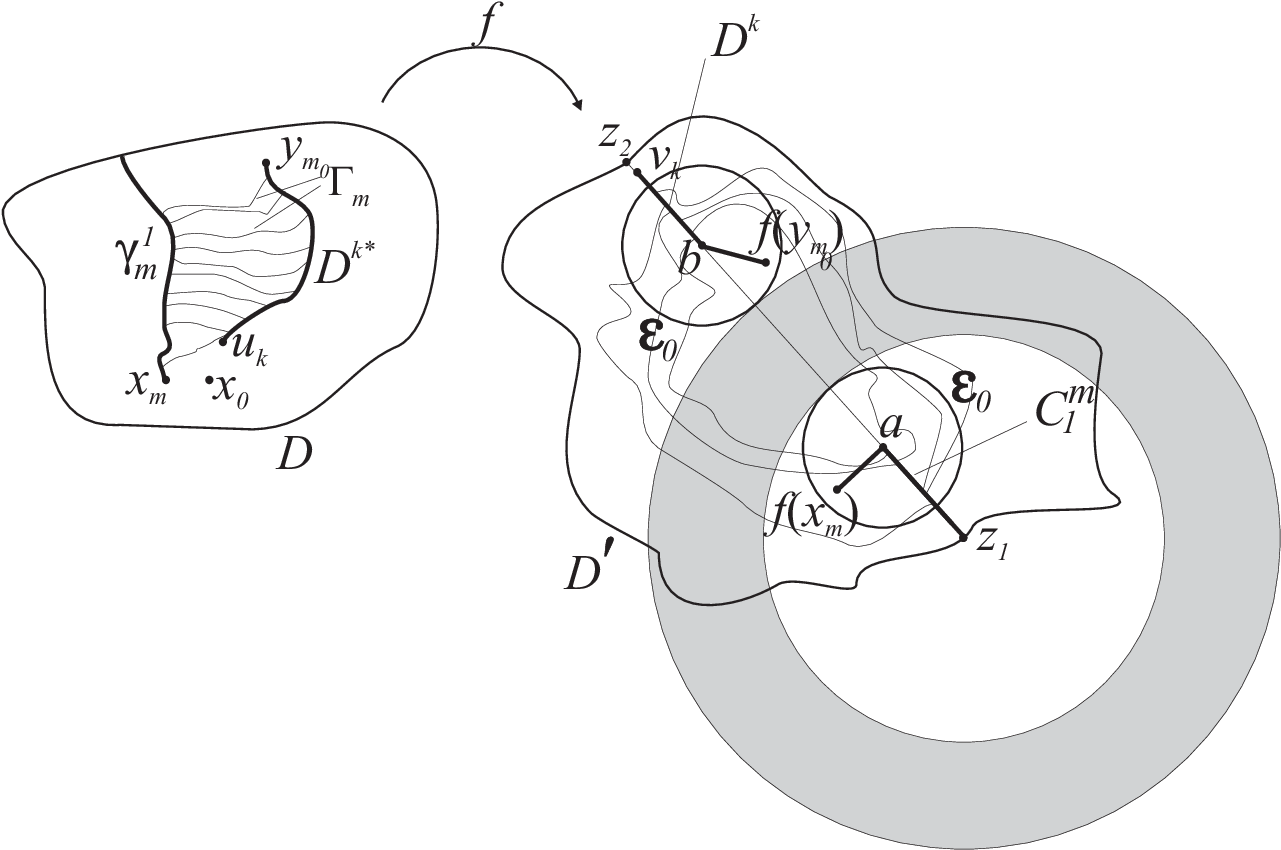}} \caption{To
the proof of Lemma~\ref{lem3}, case~2)}\label{fig5}
\end{figure}
Now we argue in the same way as in the case~1), where $\gamma_{m}^1$
and $C^{m}_1$ are constructed similarly, and $\gamma_{m}^2$ and
$C^{m}_2$ perform the mean of $D^{\,m*}$ and $D^{m},$ respectively.
Due to the definition of $\gamma_{m}^1$ and $D^{\,m*},$ there is
$\delta_0>0$ such that $d(|\gamma_{m}^1|)\geqslant \delta_0>0$ and
$d(|D^{\,m*}|)\geqslant \delta_0>0$ for any $m=1,2,\ldots .$

\medskip
Fix $P>0.$ Let $U:=B(x_0, \delta_0/2),$ and let $V$ be a
neighborhood of the same point $x_0,$ which corresponds to
Lemma~\ref{lem2} and Remark~\ref{rem1}. Reasoning similarly to the
case~1), we may obtain that
\begin{equation}\label{eq11I}
M(\Gamma(|\gamma_{m}^1|, |D^{m\,*}|, D\setminus \{x_0\}))>P\,\qquad
m\geqslant m_0
\end{equation}
for $p=n.$ If $p>n,$
\begin{equation}\label{eq1M}
M_p(\Gamma(|\gamma_{m}^1|, |D^{m\,*}|, D\setminus\{x_0\}))\geqslant
P^{\frac{p}{n}}\cdot (m(D))^{\frac{n-p}{n}}\,.
\end{equation}
Let us show that the relations~(\ref{eq11I}) and~(\ref{eq1M}) for
$p=n$ and $p>n,$ respectively, are impossible (namely, they
contradict with the definition of~$f$
in~(\ref{eq2*A})--(\ref{eqA2})). Indeed, arguing similarly to the
case~1), due to the fact that, the paths $C^{m}_1$ and $D^{m\,*}$
satisfy the conditions $|C_1^m|\subset B(z_1,
|a-z_1|+\varepsilon_0)$ and $D^{m\,*}\subset |C_2^{m_{0}}|\subset
D^{\,\prime}\setminus \overline{B(z_1, |b-z_1|-\varepsilon_0)}$
(see~(\ref{eq16}) and~(\ref{eq17})) we obtain that
\begin{equation}\label{eq14E}
M_p(\Gamma(|\gamma_{m}^1|, |D^{m\,*}|, D\setminus\{x_0\}))\leqslant
\frac{\omega_{n-1}}{I^{p-1}}<\infty\,.
\end{equation}
The relation~(\ref{eq14E}) contradicts with~(\ref{eq11I}) for $p=n$
and~(\ref{eq1M}) for $p>n$. This contradiction completes the
considerations of the case~2).

\medskip
As we say above, {\bf the case~3)} may be considered completely
similarly to the case~2) (it differs from it by the notions of paths
and corresponding points, only).

\medskip
Finally, let us to consider {\bf the case~4)}, which reduces to the
case~1) in some degree.

\medskip
Assume that, there are $k_0, m_0\in {\Bbb N}$ such that
$h(\gamma_{k_0}^1(t), x_0)\rightarrow 0$ as $t\rightarrow
c^1_{k_0}-0$ and $h(\gamma_{m_0}^2(t), x_0)\rightarrow 0$ as
$t\rightarrow c^2_{m_0}-0.$

\medskip
Since ${\gamma}^{\,0}_1(t)\rightarrow \{x_0\}$ as $t\rightarrow
c^{k_0}_2-0,$ there is a sequence $s_k\rightarrow c^{k_0}_2-0$ such
that ${\gamma}^{\,0}_1(s_k)\rightarrow x_0$ as $k\rightarrow\infty.$
Set $w_k:={\gamma}^{\,0}_1(s_k)$ and
$\omega_k:=f({\gamma}^{\,0}_1(s_k)).$ Let also
\begin{equation}\label{eq15E}
E^{\,k*}:={{\gamma}^{\,0}_1}|_{[0, s_k]}\,, \quad
E^k:={C^{k_0\,}_1}|_{[0, s_k]}\,,\quad k=1,2,\ldots\,.
\end{equation}
Similarly, since $\gamma_{m_0}^2(t)\rightarrow x_0$ as $t\rightarrow
c^{m_0}_2-0,$ there exists a sequence $t_k\rightarrow c^{m_0}_2-0$
such that $\gamma_{m_0}^2(t_k)\rightarrow x_0.$ Set
$u_k:=\gamma_{m_0}^2(t_k)$ і $v_k:=f(\gamma_{m_0}^2(t_k)).$ Let
\begin{equation}\label{eq15C}
D^{\,k*}:={\gamma_{m_0}^2}|_{[0, t_k]}\,, \quad
D^k:={C^{m_0\,}_2}|_{[0, t_k]}\,,\quad k=1,2,\ldots\,,
\end{equation}
see Figure~\ref{fig6}.
\begin{figure}[h]
\centerline{\includegraphics[scale=0.5]{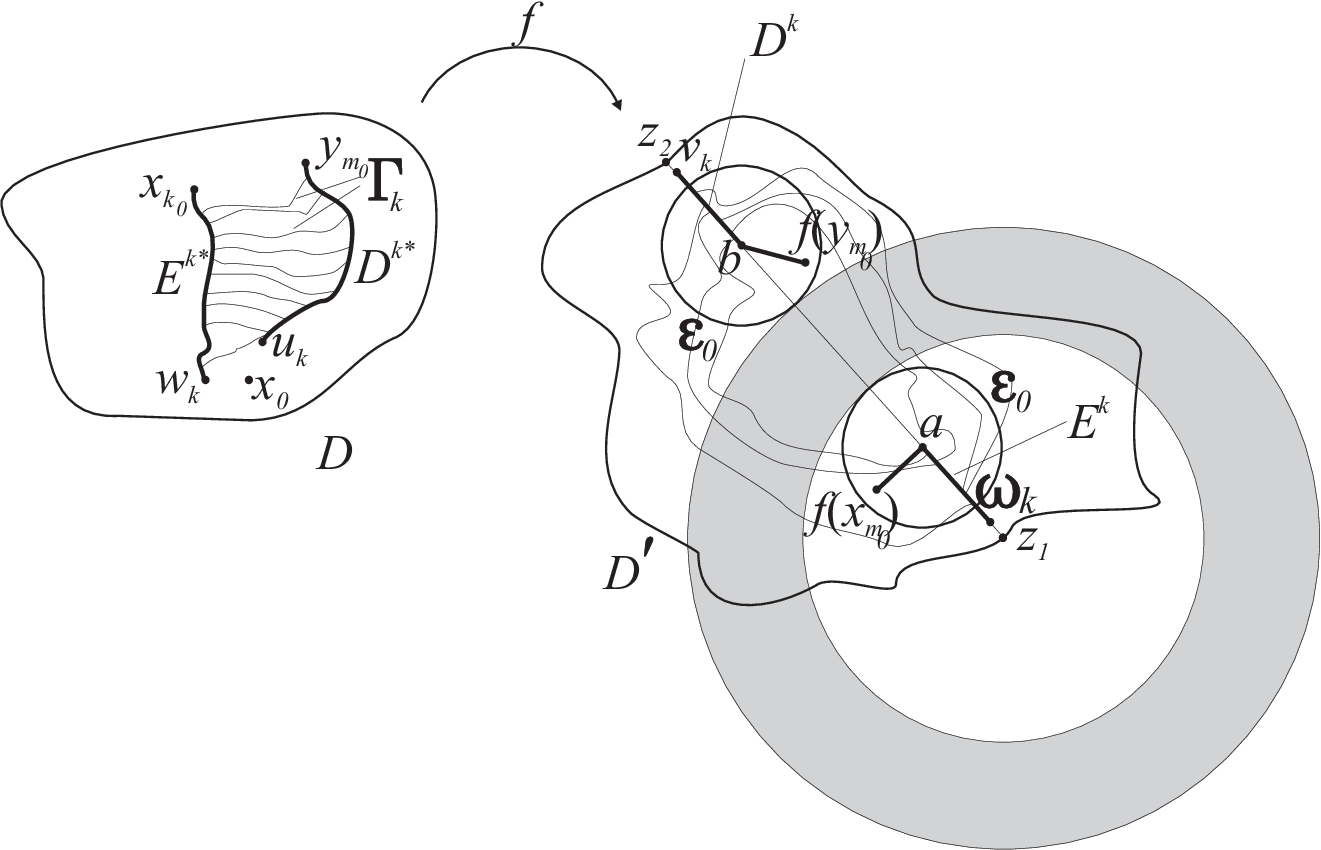}} \caption{To
the proof of Lemma~\ref{lem1}, case~4)}\label{fig6}
\end{figure}
Now we are reasoning similarly to the case~1), where we consider
$\gamma_m^1$ and $C^{m}_1$ instead of $E^{\,m*}$ and $E^{m},$
respectively, and $\gamma_m^2$ and $C^{m}_2$ instead of $D^{\,m*}$
and $D^{m},$ respectively. Due to the definition of $E^{\,k*}$ and
$D^{\,k*},$ there is $\delta_0>0$ such that $d(|E^{\,k*}|)\geqslant
\delta_0>0$ and $d(|D^{\,k*}|)\geqslant \delta_0>0$ for any $k\in
{\Bbb N}.$

\medskip
Fix $P>0.$ Let $U:=B(x_0, \delta_0/2),$  and let $V$ be a
neighborhood of the same point $x_0,$ which corresponds to
Lemma~\ref{lem2} and Remark~\ref{rem1}. Arguing similarly to the
case~1), we obtain that
\begin{equation}\label{eq11J}
M(\Gamma(|E^{k\,*}|, |D^{k\,*}|, D\setminus \{x_0\}))>P\,,\quad
k\geqslant K_0
\end{equation}
for $p=n.$ If $p>n,$
\begin{equation}\label{eq1N}
M_p(\Gamma(|E^{k\,*}|, |D^{k\,*}|, D\setminus\{x_0\}))\geqslant
P^{\frac{p}{n}}\cdot (m(D))^{\frac{n-p}{n}}\,.
\end{equation}
Let us to show that, the relations~(\ref{eq11J}) and~(\ref{eq1N})
for $p=n$ and $p>n,$ respectively, are impossible because they
contradict with the definition of the mapping~$f$
in~(\ref{eq2*A})--(\ref{eqA2})). Indeed, arguing similarly to the
case~1) and taking into account that, the paths $E^{\,k}$ belong
$D^{\,k}$ to disjoint balls $B(z_1, \varepsilon_0)$ and $B(z_2,
\varepsilon_0),$ we obtain that
\begin{equation}\label{eq14F}
M_p(\Gamma(\Gamma(|E^{k\,*}|, |D^{k\,*}|,
D\setminus\{x_0\}))\leqslant \frac{\omega_{n-1}}{I^{p-1}}<\infty\,.
\end{equation}
The relation~(\ref{eq14F}) contradict with~(\ref{eq11J}) for $p=n$
and~(\ref{eq1N}) for $p>n$. The contradiction obtained above
completes the consideration of the case~4). Lemma is completely
proved.~$\Box$
\end{proof}

\medskip
Lemma~\ref{lem3} admits some strengthening if $Q\in
L^1(D^{\,\prime}):=f(D\setminus\{x_0\}).$ The following is true.

\medskip
\begin{lemma}\label{lem4}
{\sl\, Let $D$ be a domain in ${\Bbb R}^n,$ $n\geqslant 2,$ let
$p\geqslant n,$ $x_0\in D$ and let $f:D\setminus\{x_0\}\rightarrow
D^{\,\prime},$ $D^{\,\prime}:=f(D\setminus\{x_0\}),$ be an open
discrete mapping that satisfies
conditions~(\ref{eq2*A})--(\ref{eqA2}) at each point $y_0\in
D^{\,\prime}.$

\medskip
Assume that $Q\in L^1(D^{\,\prime}):=f(D\setminus\{x_0\}).$ Then the
set $C(x_0, f)\cap D^{\,\prime}$ cannot contain more than one
points.}
\end{lemma}

\medskip
\begin{proof}
The proof of this statement is very similar to the proof of the
previous lemma, however, for the sake of completeness, we have given
it in its entirety. Suppose the opposite. Then there exist at least
two sequences $x_m, y_m\in D\setminus\{x_0\},$ $x_m, y_m\rightarrow
x_0$ as $m\rightarrow\infty$ such that $f(x_m)\rightarrow a$ and
$f(y_m)\rightarrow b$ as $m\rightarrow\infty,$ where $a, b\in
D^{\,\prime}$ and $a\ne b.$ We may assume that $f(x_m)\in B(a,
\varepsilon_0)\subset D$ and $f(y_m)\in B(b, \varepsilon_0)\subset
D$ for all $m\in {\Bbb N},$ where $\varepsilon_0>0$ is such that
\begin{equation}\label{ewq13***A}
\overline{B(a, \varepsilon_0)}\cap \overline{B(b,
\varepsilon_0)}=\varnothing\,.
\end{equation}
In particular, it follows from ~(\ref{eq13***}) that
\begin{equation}\label{eq13D}
\varepsilon_0<\frac{|a-b|}{2}\,.
\end{equation}
Consider the line $$r=r(t)=a+(b-a)t\,,\quad-\infty<t<\infty\,,$$
which passes through points $a$ and $b,$ see Figure~\ref{fig4A}.
\underline{}\begin{figure}[h]
\centerline{\includegraphics[scale=0.5]{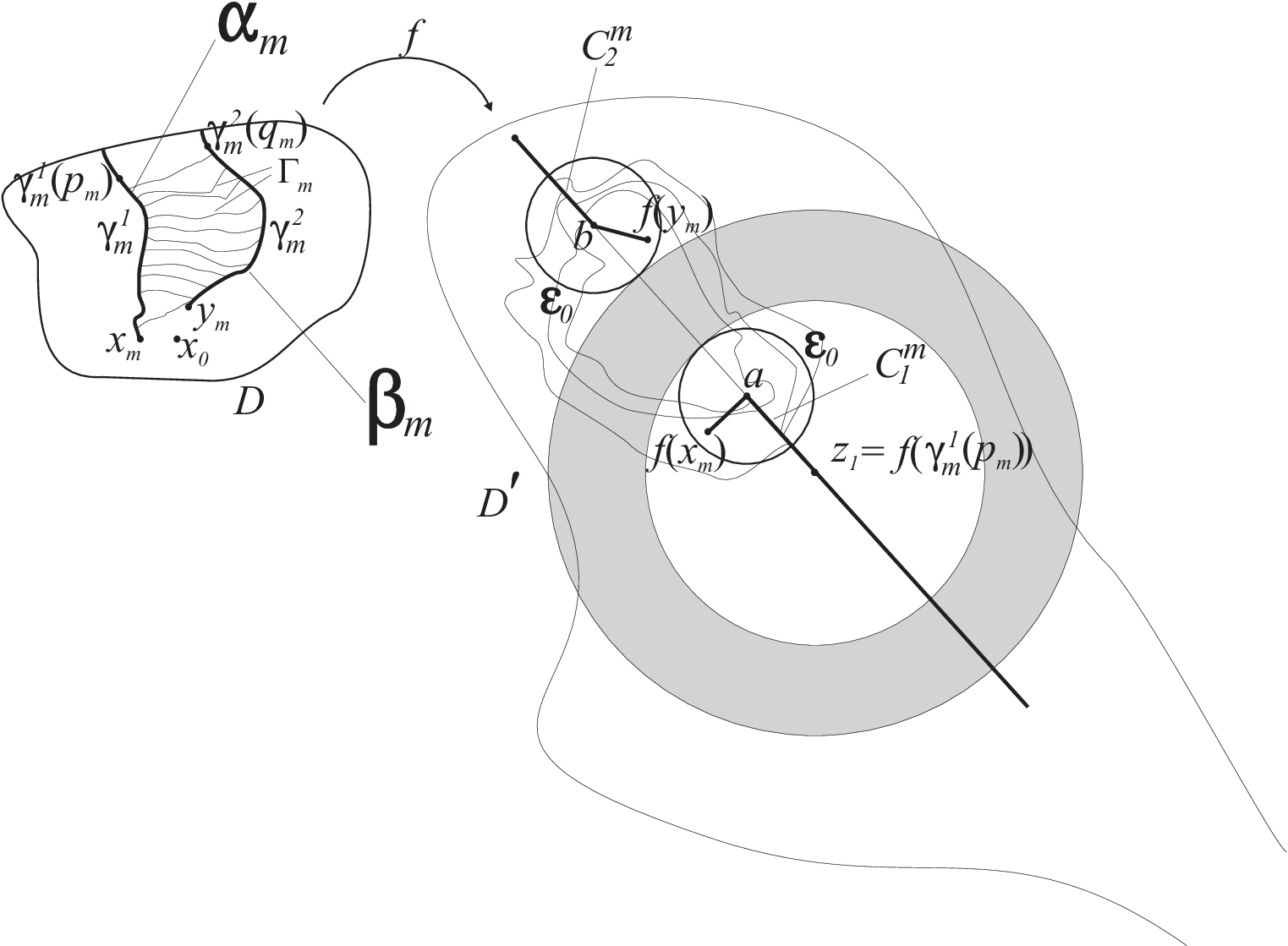}} \caption{To
the proof of Lemma~\ref{lem4}, case~1) }\label{fig4A}
\end{figure}
Let us denote by $C^m_1$ the broken line $\overline{f(x_m)a\infty},$
which is a union of the segment $\overline{f(x_m)a}$ with a ray
$\{r=r(t), t\leqslant 0\},$ and by $C^m_2$ the broken line
$\overline{f(y_m)b\infty},$ which is a union of the segment
$\overline{f(y_m)b}$ with the ray $\{r=r(t), t\geqslant 1\}.$ We
parameterize these broken lines in the following way:
$$C^m_1:[0, 1)\rightarrow {\Bbb R}^n, \quad C^m_1(0)=f(x_m)\,,$$
$$C^m_2:[0, 1)\rightarrow {\Bbb R}^n, \quad C^m_2(0)=f(y_m)\,.$$
Let $\gamma_m^1:[0, c^1_m)\rightarrow D\setminus\{x_0\},$
$0<c^1_m\leqslant 1,$ be a maximal $f$-lifting of the path $C_1^{m}$
with the origin at the point $x_m,$ which exists by
Proposition~\ref{pr3}. Let us to prove that, the situation when
$\gamma_m^1\rightarrow \omega_1\in D\setminus\{x_0\}$ as
$t\rightarrow c^1_m,$ is not possible. Indeed, by
Proposition~\ref{pr3} we would have that $c^1_m=1$ and
$f(\omega_1)=\lim\limits_{t\rightarrow 1-0}C_1^m(t).$ Then, from one
side, $f(\omega_1)\in D^{\,\prime}$ by the openness of the mapping
$f,$ and on the other hand, $f(\omega_1)=z_1$ by definition of
$C_1^m(t).$ Since $z_1\not\in D^{\,\prime},$ we obtain a
contradiction. So, the case $\gamma_m^1(t)\rightarrow \omega_1\in
D\setminus\{x_0\}$ as $t\rightarrow c^1_m$ is impossible. Again by
Proposition~\ref{pr3}
\begin{equation}\label{eq3A}
h(\gamma_m^1(t), \partial (D\setminus\{x_0\}))\rightarrow 0
\end{equation}
as $t\rightarrow c^1_m-0.$ Similarly, let $\gamma_m^2:[0,
c^2_m)\rightarrow D\setminus\{x_0\},$ $0<c^2_m\leqslant 1,$ be a
maximal $f$-lifting of the path $C_2^m$ starting at the point $y_m,$
which exists by Proposition~\ref{pr3}. Similarly to~(\ref{eq3A}), we
obtain that
\begin{equation}\label{eq4A}
h(\gamma_m^2(t), \partial D\setminus\{x_0\})\rightarrow 0
\end{equation}
as $t\rightarrow c^2_m-0.$

\medskip
Let us establish that, for any point $z_1=r(t),$ $t<0,$
\begin{equation}\label{eq16A}
\overline{f(x_m)az_1}\subset B(z_1, |a-z_1|+\varepsilon_0)\,,
\end{equation}
where $\overline{f(x_m)az_1}$ denotes the broken line that sequently
joins the points $f(x_m),$ $a$ and $z_1.$

\medskip Indeed, if $x\in \overline{f(x_m)az_1},$ then either $x\in \overline{az_1}$
(where $\overline{az_1}$ is the segment joining the points $a$ and
$z_1$), or $x\in \overline{f(x_m)a}\subset B(a, \varepsilon_0).$ If
$x\in \overline{az_1},$ then, obviously, $|x-z_1|\leqslant
|a-z_1|<|a-z_1|+\varepsilon_0;$ therefore, $x\in B(z_1,
|a-z_1|+\varepsilon_0|).$ Now, if $x\in \overline{f(x_m)a}\subset
B(a, \varepsilon_0),$ then
$|x-a|<\varepsilon_0<|a-z_1|+\varepsilon_0|.$ Therefore also $x\in
B(z_1, |a-z_1|+\varepsilon_0|).$ The inclusion~(\ref{eq16A}) is
established.

\medskip
Let us now establish that, for any $m\in {\Bbb N}$ and for any
points $z_1=r(t),$ $t<0,$
\begin{equation}\label{eq17A}
|C_2^m|\subset {\Bbb R}^n\setminus \overline{B(z_1,
|b-z_1|-\varepsilon_0)}\,.
\end{equation}
Let $x\in |C_2^m|.$ Then either $x\in \overline{bz_2}$ (where
$\overline{bz_2}$ is a segment joining the points $b$ and $z_2$), or
$x\in \overline{f(y_m)b}\subset B(b, \varepsilon_0).$

Assume that $x\in \overline{bz_2}.$ By the construction
$x=a+(b-a)t_1,$ $t_1\geqslant 1,$ and $z_1=a+(b-a)t_2,$
$t_2\leqslant 0.$ Then
$$|x-z_1|=|a+(b-a)t_1-a-(b-a)t_2|=|(b-a)(t_1-t_2)|\geqslant$$
$$\geqslant |b-a|(1-t_2)= |z_1-b|\geqslant
|b-z_1|-\varepsilon_0\,,$$
therefore $x\in D^{\,\prime}\setminus \overline{B(z_1,
|b-z_1|-\varepsilon_0|)}.$

Now, let $x\in \overline{f(y_m)b}\subset B(b, \varepsilon_0).$ Then,
by the triangle inequality,
$$|x-z_1|\geqslant |z_1-b|-|b-x|\geqslant |b-z_1|-\varepsilon_0\,.$$
The relation~(\ref{eq17A}) is also established.

\medskip
It follows from~(\ref{eq3A}) and~(\ref{eq4A}) that, the following
four cases are possible:

\medskip
1) for any $m\in {\Bbb N},$ $h(\gamma_{m}^1(t), \partial
D)\rightarrow 0$ as $t\rightarrow c^1_{m}-0$ and $h(\gamma_{m}^2(t),
\partial D)\rightarrow
0$ as $t\rightarrow c^2_{m}-0;$

\medskip
2) for any $m\in {\Bbb N},$ $h(\gamma_{m}^1(t), \partial
D)\rightarrow 0$ as $t\rightarrow c^1_{m}-0,$ but there exists
$m_0\in {\Bbb N}$ such that $h(\gamma_{m_0}^2(t), x_0)\rightarrow 0$
as $t\rightarrow c^2_{m_0}-0;$

\medskip
3) there exists $k_0\in {\Bbb N}$ such that $h(\gamma_{k_0}^1(t),
x_0)\rightarrow 0$ as $t\rightarrow c^1_{k_0}-0,$ in addition, for
any $m\in {\Bbb N},$ $h(\gamma_{m}^2(t),
\partial D)\rightarrow 0$ as $t\rightarrow c^2_{m}-0;$

\medskip
4) there are numbers $k_0, m_0\in {\Bbb N}$ such that
$h(\gamma_{k_0}^1(t), x_0)\rightarrow 0$ as $t\rightarrow
c^1_{k_0}-0$ and $h(\gamma_{m_0}^2(t), x_0)\rightarrow 0$ as
$t\rightarrow c^2_{m_0}-0.$

\medskip
Let us consider each case separately. Note that case~2) is some
''symmetric'' to case~3), so that these two cases are sufficient
consider only one, for example, case~2).

\medskip
Consider the {\bf case~1)}: for any $m\in {\Bbb N},$
$h(\gamma_{m}^1(t), \partial D)\rightarrow 0$ as $t\rightarrow
c^1_{m}-0$ and $h(\gamma_{m}^2(t),
\partial D)\rightarrow
0$ as $t\rightarrow c^2_{m}-0$ (see Figure~\ref{fig4A}). Then there
are sequences $p_m\rightarrow c^1_{m}-0$ and $q_m\rightarrow
c^2_{m}-0$ as $m\rightarrow\infty$ such that $h(\gamma_{m}^1(p_m),
\partial D)\rightarrow 0$ as $m\rightarrow \infty$ and
$h(\gamma_{m}^2(q_m),
\partial D) \rightarrow
0$ at $m\rightarrow \infty.$ Set
\begin{equation}\label{eq21A}
\alpha_m:=\gamma_{m}^1|_{[0, p_m]}\,,\qquad
\beta_m:=\gamma_{m}^2|_{[0, q_m]}\,. \end{equation}
Then there exists $\delta_0>0$ such that $d(|\alpha_m|)\geqslant
\delta_0>0$ and $d(|\beta_m|)\geqslant \delta_0>0$ for all $m\in
{\Bbb N}.$ Let us fix the number $P>0.$ Let $U:=B(x_0, \delta_0/2),$
and let $V$ be a neighborhood of the same point $x_0,$ which
corresponds to Lemma~\ref{lem2} and Remark~\ref{rem1}. Since by the
assumption $x_m, y_m\in D\setminus \{x_0\},$ $m=1,2,\ldots,$ then we
may find the number $m_0\in{\Bbb N}$ such that $x_m, y_m\in V$ for
all $m\geqslant m_0.$
Note that, for $m\geqslant m_0,$
\begin{equation}\label{eq10C}
|\alpha_m|\cap \partial V\ne\varnothing, \quad |\beta_m|\cap
\partial V\ne\varnothing\,.
\end{equation}
Indeed, $x_m\in |\alpha_m|,$ $y_m\in |\beta_m|,$ thus
$|\alpha_m|\cap V\ne\varnothing\ne |\beta_m|\cap V$ at $m\geqslant
m_0.$ In addition, ${\rm diam}\,V\leqslant {\rm diam}\,U=r_0$ and,
since $d(|\alpha_m|)\geqslant \delta_0>0$ and $d(|\beta_m|)\geqslant
\delta_0>0$ for all $m\in {\Bbb N},$ by Proposition~\ref{pr2} we
obtain the relation~(\ref{eq10C}). Similarly, we may prove that
\begin{equation}\label{eq11K}
|\alpha_m|\cap \partial U\ne\varnothing, \quad |\beta_m|\cap
\partial U\ne\varnothing\,.
\end{equation}
Then,by Lemma~\ref{lem2} and Remark~\ref{rem1}
\begin{equation}\label{eq11L}
M(\Gamma(|\alpha_m|, |\beta_m|, D\setminus \{x_0\}))>P\,,\qquad
m\geqslant m_0.
\end{equation}
If $p>n,$ then by H\"{o}lder inequality, for every $\rho\in {\rm
adm}\,\Gamma(|\alpha_m|, |\beta_m|, D\setminus\{x_0\})$ we obtain
that
\begin{equation}\label{eq1O}
M(\Gamma(|\alpha_m|, |\beta_m|, D\setminus\{x_0\}))\leqslant
\int\limits_{D}\rho^n(x)\,dm(x)\leqslant
\left(\int\limits_{D}\rho^p(x)\,dm(x)\right)^{\frac{n}{p}}\cdot
(m(D))^{\frac{p-n}{p}}\,.
\end{equation}
Letting in~(\ref{eq1O}) to $\inf$ over all $\rho\in {\rm
adm}\,\Gamma(|\alpha_m|, |\beta_m|, D\setminus\{x_0\}),$ we obtain
that
$$M(\Gamma(|\alpha_m|, |\beta_m|,
D\setminus\{x_0\}))\leqslant$$
\begin{equation}\label{eq1P} \leqslant
\int\limits_{D}\rho^n(x)\,dm(x)\leqslant
\left(M_p(\Gamma(|\alpha_m|, |\beta_m|,
D\setminus\{x_0\}))\right)^{\frac{n}{p}}\cdot
(m(D))^{\frac{p-n}{n}}\,.
\end{equation}
By~(\ref{eq11L}) and~(\ref{eq1P}) it follows that
\begin{equation}\label{eq1R}
M_p(\Gamma(|\alpha_m|, |\beta_m|, D\setminus\{x_0\}))\geqslant
P^{\frac{p}{n}}\cdot (m(D))^{\frac{n-p}{n}}\,.
\end{equation}
We will show that, the relations~(\ref{eq11L}) and~(\ref{eq1R}) for
$p=n$ and $p>n,$, respectively, are impossible (in particular, each
of them contradicts the definition of the mapping~$f$
in~(\ref{eq2*A})-(\ref{eqA2})).

\medskip
Let $\Gamma_*=\Gamma(\overline{f(x_m)af(\gamma^1_m(p_m))}, |C^m_2|,
D^{\,\prime}),$ where $\overline{f(x_m)af(\gamma^1_m(p_m))}$ denotes
a broken line that sequently joins the points $f(x_m),$ $a$ and
$f(\gamma^1_m(p_m)).$ Observe that
\begin{equation}\label{eq3F}
\Gamma_*>\Gamma(S(z_1, r_1|), S(z_1, r_2), A(z_1, r_1, r_2)\,,
\end{equation}
where $$z_1=f(\gamma^1_m(p_m))\,,$$
\begin{equation}\label{eq20A}
r_1=|a-z_1|+\varepsilon_0\,, \qquad r_2=|b-z_1|-\varepsilon_0\,.
\end{equation}
First of all, let us clarify that $r_2>r_1.$ Indeed, if
$z_1=a+(b-a)t_2$ for of some $t_2\leqslant 0,$ then
$$r_2-r_1=|b-z_1|-\varepsilon_0-|a-z_1|-\varepsilon_0=$$
\begin{equation}\label{eq21}
=|b-a-(b-a)t_2|-|a-a-(b-a)t_2|-2\varepsilon_0=
\end{equation}
$$=|b-a|(1-t_2)+|b-a|t_2-2\varepsilon_0=|b-a|-
2\varepsilon_0\geqslant 0\,,$$
see~(\ref{eq13D}).

Now, let $\gamma\in \Gamma_*,$ $\gamma:[0, 1]\rightarrow {\Bbb
R}^n.$ Due to~(\ref{eq16A}) and~(\ref{eq17A}), $\gamma(0)\in
\overline{f(x_m)af(\gamma^1_m(p_m))}\subset B(z_1, r_1)$ and
$\gamma(1)\in |C^m_2|\subset D^{\,\prime}\setminus \overline{B(z_1,
r_2)}.$ Now, by Proposition~\ref{pr2} there is $t_1\in (0, 1)$ such
that $\gamma(t_1)\in S(z_1, r_1).$ Without loss of restriction, it
is possible assume that $|\gamma(t)-z_1|>r_1$ when $t>t_1.$ Next,
since $\gamma(t_1)\in S(z_1, r_1)$ and $\gamma(1)\in |C^m_2|\subset
{\Bbb R}^n\setminus B(z_1, r_2),$ by Proposition~\ref{pr2} there is
$t_2\in (t_1, 1)$ such that $\gamma(t_2)\in S(z_1, r_2).$ Without
loss of generalization, we may assume that $|\gamma(t)-z_1|<r_2$ at
$t_1<t<t_2.$ Therefore, $\gamma|_{[t_1, t_2]}$ is a subpath of
$\gamma$ which belongs to~$\Gamma(S(z_1, r_1), S(z_1, r_2), A(z_1,
r_1, r_2)).$ The relation~(\ref{eq3F}) is proved.

\medskip
Let us establish now that
\begin{equation}\label{eq5D}
\Gamma(|\alpha_m|, |\beta_m|, D\setminus\{x_0\})>\Gamma_f(z_1, r_1,
r_2)\,,
\end{equation}
where, as above, $r_1$ and $r_2$ are defined in~(\ref{eq20A}).
Indeed, if $\gamma:[0, 1]\rightarrow D\setminus\{x_0\}$ belongs to
$\Gamma(|\alpha_m|, |\gamma^2_m|, D\setminus\{x_0\}),$ then
$f(\gamma)$ belongs to $D^{\,\prime},$ and $f(\gamma(0))\in
\overline{f(x_m)af(\gamma^1_m(p_m))}$ and $f(\gamma(1))\in
|C_2^{m}|,$ i.e., $f(\gamma)\in \Gamma_*.$ Then, taking into account
mentioned above, by~(\ref{eq3F}) $f(\gamma)$ has a subpath
$f(\gamma)^{\,*}:=f(\gamma)|_{[t_1, t_2]},$ $0\leqslant
t_1<t_2\leqslant 1,$ which belongs to the family $\Gamma(S(z_1,
r_2), S(z_1, r_1), A(z_1, r_1, r_2)).$ Then
$\gamma^*:=\gamma|_{[t_1, t_2]}$ is subpath of $\gamma$ which
belongs to~$\Gamma_f(z_1, r_1, r_2),$ as required.

\medskip
In turn, by~(\ref{eq5D}) we have that
$$M_p(\Gamma(|\alpha_m|, |\beta_m|, D\setminus\{x_0\}))\leqslant$$
\begin{equation}\label{eq11M}
\leqslant M_p(\Gamma_f(z_1, r_1, r_2))\leqslant \int\limits_{A}
Q(y)\cdot \eta^p (|y-z_1|)\, dm(y)\,,
\end{equation}
where $A=A(z_1, r_1, r_2)$ and $\eta$ is an arbitrary Lebesgue
measurable function satisfying the relation~(\ref{eqA2}).
Now let us put
$$\eta(t)= \left\{
\begin{array}{rr}
\frac{1}{r_2-r_1}, & t\in [r_1, r_2],\\
0, & t\not\in [r_1, r_2]\,.
\end{array}
\right. $$
Note that, $\eta$ satisfies the relation~(\ref{eqA2}). Then
\begin{equation}\label{eq14***A}
M_p(\Gamma(|\alpha_m|, |\beta_m|, D\setminus\{x_0\}))\leqslant
M_p(\Gamma_{f}(z^1, r_1, r_2))\leqslant
\frac{1}{(r_2-r_1)^p}\int\limits_{D^{\,\prime}}
Q(y)\,dm(y)<\infty\,.
\end{equation}
The relation~(\ref{eq14***A}) contradicts with~(\ref{eq11L}) for
$p=n$ and~(\ref{eq1R}) for $p>n$. The resulting contradiction
completes the consideration of the case~1).

\medskip
Let us consider the {\bf case 2):} for any $m\in {\Bbb N},$
$h(\gamma_{m}^1(t), \partial D)\rightarrow 0$ as $t\rightarrow
c^1_{m}-0,$ but there exists $m_0\in {\Bbb N}$ such that
$h(\gamma_{m_0}^2(t), x_0)\rightarrow 0$ as $t\rightarrow
c^2_{m_0}-0;$ see Figure~\ref{fig5A}.
\begin{figure}[h]
\centerline{\includegraphics[scale=0.5]{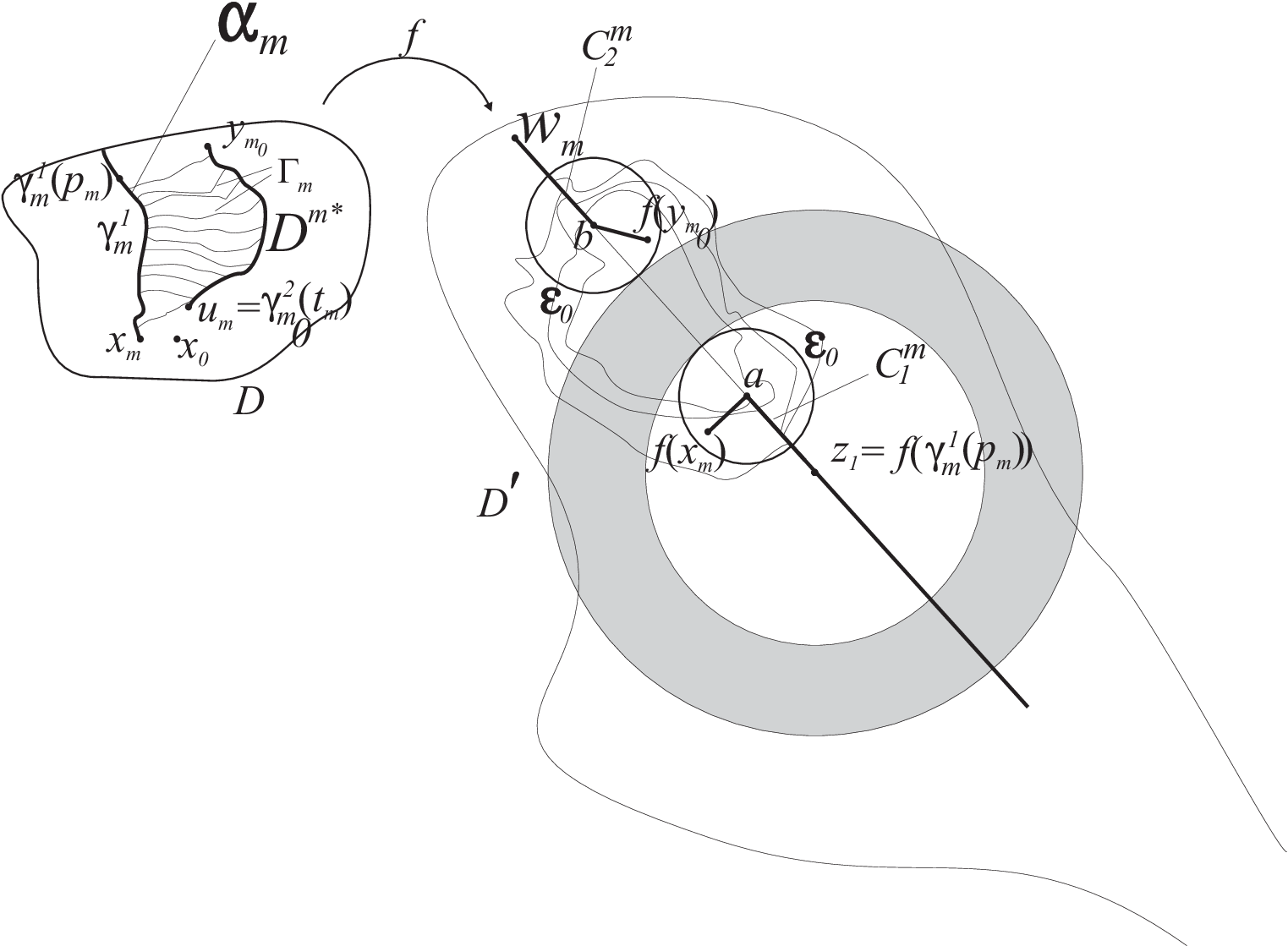}} \caption{To
the proof of Lemma~\ref{lem4}, case~2)}\label{fig5A}
\end{figure}
Since $\gamma_{m_0}^2(t)\rightarrow \{x_0\}$ as $t\rightarrow
c^{m_0}_2-0,$ there is a sequence $t_k\rightarrow c^{m_0}_2-0$ such
that $\gamma_{m_0}^2(t_k)\rightarrow x_0$ as $k\rightarrow\infty.$
We put $u_k:=\gamma_{m_0}^2(t_k)$ and $v_k:=f(\gamma_{m_0}^2(t_k)).$
Let also
\begin{equation}\label{eq15F}
D^{\,k*}:={\gamma_{m_0}^2}|_{[0, t_k]}\,, \quad
D^k:=f\left({\gamma_{m_0}^2}|_{[0, t_k]}\right)\,,\quad
k=1,2,\ldots\,.
\end{equation}
In addition, let a sequence of points $p_m,$ $m=1,2,\ldots,$ and a
path $\alpha_m,$ $m=1,2,\ldots, $ are defined in the same way as in
case~1).

By the definitions of $\alpha_m$ and $D^{\,m*},$ there exists
$\delta_0>0$ such that $d(\alpha_m)\geqslant \delta_0>0$ and
$d(|D^{\,m*}|)\geqslant \delta_0>0$ for all $m=1,2,\ldots .$

\medskip
Let us fix $P>0.$ Let $U:=B(x_0, \delta_0/2),$ and let $V$ be a
neighborhood of the same point $x_0,$ which corresponds to
Lemma~\ref{lem2} and Remark~\ref{rem1}. Reasoning in the same way as
in case~1), we obtain that
\begin{equation}\label{eq11N}
M(\Gamma(|\alpha_m|, |D^{m\,*}|, D\setminus \{x_0\}))>P\,,\qquad
m\geqslant m_0
\end{equation}
for $p=n.$ If $p>n,$
\begin{equation}\label{eq1S}
M_p(\Gamma(|\alpha_m|, |D^{m\,*}|, D\setminus\{x_0\}))\geqslant
P^{\frac{p}{n}}\cdot (m(D))^{\frac{n-p}{n}}\,.
\end{equation}
Let us show that, the relations~(\ref{eq11N}) and~(\ref{eq1S}) for
$p=n$ and $p>n,$ respectively, are impossible (in particular, each
of them contradicts definition of mapping~$f$
in~(\ref{eq2*A})--(\ref{eqA2})). Really, reasoning in the same way
as in case~1), due to the relations~(\ref{eq16A}) and~(\ref{eq17A})
we have that $f(\alpha_m)\subset B(z_1, r_1)$ and
$|D^{m\,}|=|f(D^{m\,*})|\subset {\Bbb R}^n\setminus \overline{B(z_1,
r_2)},$ where
$$z_1=f(\gamma^1_m(p_m))\,,$$
\begin{equation}\label{eq20B}
r_1=|a-z_1|+\varepsilon_0\,, \qquad r_2=|b-z_1|-\varepsilon_0\,.
\end{equation}
Similarly to~(\ref{eq5D}), we have:
\begin{equation}\label{eq5E}
\Gamma(|\alpha_m|, D^{m\,*}, D\setminus\{x_0\})>\Gamma_f(z_1, r_1,
r_2)\,.
\end{equation}
In turn, by~(\ref{eq5E}) we have the following:
$$M_p(\Gamma(|\alpha_m|, |\beta_m|, D\setminus\{x_0\}))\leqslant$$
\begin{equation}\label{eq11P}
\leqslant M_p(\Gamma_f(z_1, r_1, r_2))\leqslant \int\limits_{A}
Q(y)\cdot \eta^p (|y-z_1|)\, dm(y)\,,
\end{equation}
where $A=A(z_1, r_1, r_2)$ and $\eta$ is an arbitrary Lebesgue
measurable function satisfying the relation~(\ref{eqA2}).
Now let us put
$$\eta(t)= \left\{
\begin{array}{rr}
\frac{1}{r_2-r_1}, & t\in [r_1, r_2],\\
0, & t\not\in [r_1, r_2]\,.
\end{array}
\right. $$
Note that, $\eta$ satisfies the relation~(\ref{eqA2}). Then
\begin{equation}\label{eq14****A}
M_p(\Gamma)\leqslant M_p(\Gamma_{f}(z^1, r_1, r_2))\leqslant
\frac{1}{(r_2-r_1)^p}\int\limits_{D^{\,\prime}}
Q(y)\,dm(y)<\infty\,.
\end{equation}
The relation~(\ref{eq14****A}) contradicts~(\ref{eq11N}) for $p=n$
and~(\ref{eq1S}) for $p>n$. The resulting contradiction completes
the consideration of the case~2).

\medskip
{\bf Case 3):} there exists $k_0\in {\Bbb N}$ such that
$h(\gamma_{k_0}^1(t), x_0)\rightarrow 0$ at $t\rightarrow
c^1_{k_0}-0,$ in addition, if any $m\in {\Bbb N},$
$h(\gamma_{m}^2(t),
\partial D)\rightarrow 0$ as $t\rightarrow c^2_{m}-0,$
see Figure~\ref{fig5B}.
\begin{figure}[h]
\centerline{\includegraphics[scale=0.5]{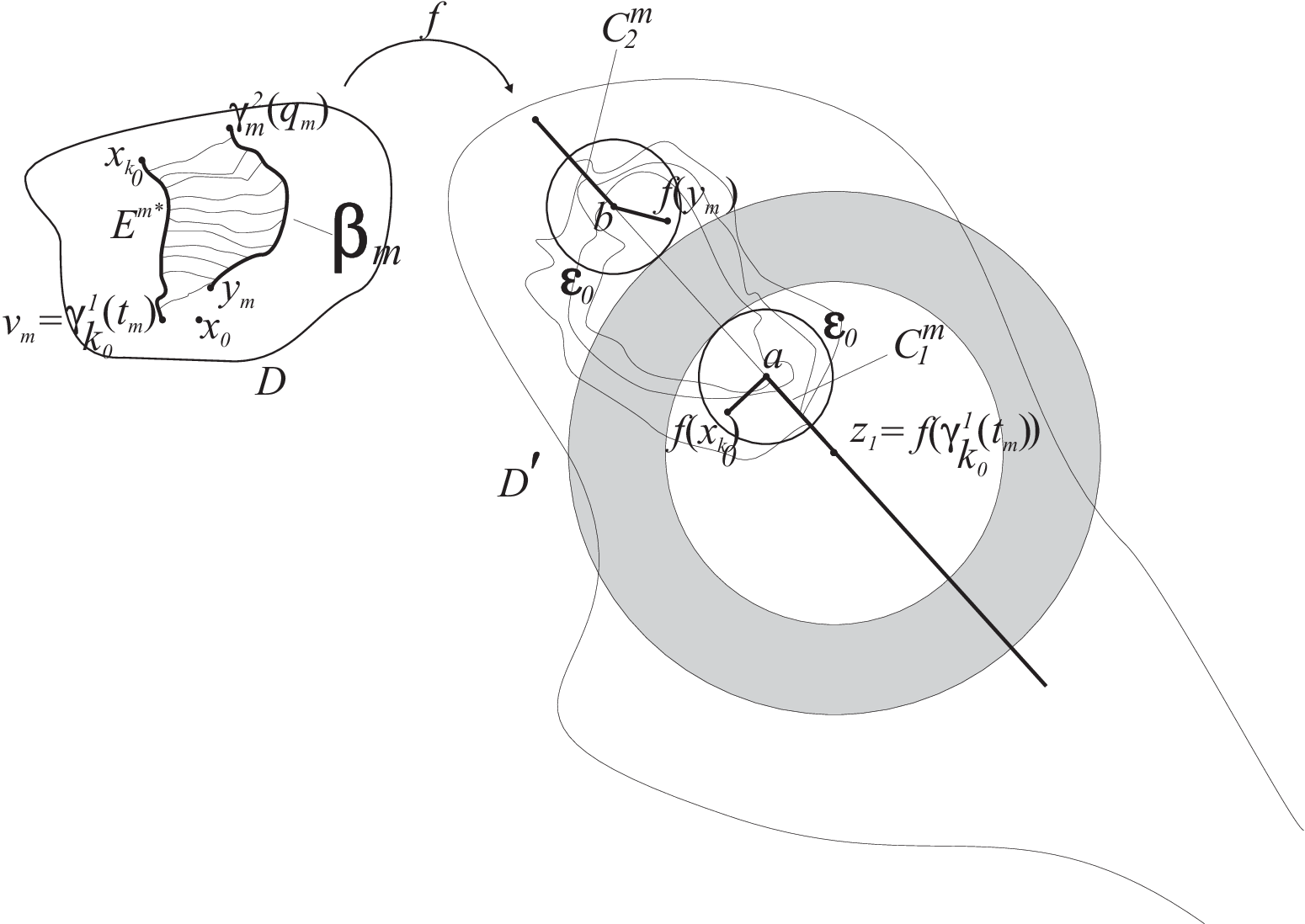}} \caption{To
the proof of Lemma~\ref{lem4}, case~3)}\label{fig5B}
\end{figure}
Since $\gamma_{k_0}^1(t)\rightarrow \{x_0\}$ as $t\rightarrow
c^{k_0}_1-0,$ there is a sequence $t_k\rightarrow c^{k_0}_1-0,$
$k\rightarrow\infty,$ such that $\gamma_{k_0}^1(t_k)\rightarrow x_0$
as $k\rightarrow\infty.$ Put $v_k:=\gamma_{m_0}^1(t_k)$ and
$z_1=z_1(k):=f(\gamma_{m_0}^1(t_k)).$ Also let
\begin{equation}\label{eq15D}
E^{\,k*}:={\gamma_{k_0}^1}|_{[0, t_k]}\,, \quad
E^k:=f\left({\gamma_{k_0}^1}|_{[0, t_k]}\right)\,,\quad
k=1,2,\ldots\,,
\end{equation}
In addition, let the sequence of points $q_m,$ $m=1,2,\ldots, $ and
the path $\beta_m,$ $m=1,2,\ldots, $ are defined in the same way as
in case~1).

By the definitions of $E^{\,m*}$ and $\beta_m$, there exists
$\delta_0>0$ such that $d(|E^{\,m*}|)\geqslant \delta_0>0$ and
$d(|\beta_m|)\geqslant \delta_0>0$ for all $m=1,2,\ldots .$

\medskip
Let us fix $P>0.$ Let $U:=B(x_0, \delta_0/2),$ and let $V$ be a
neighborhood of the same point $x_0,$ which corresponds to
Lemma~\ref{lem2} and Remark~\ref{rem1}. Reasoning in the same way as
in case~1), we have that
\begin{equation}\label{eq11O}
M(\Gamma(|E^{m\,*}|, |\beta_m|, D\setminus \{x_0\}))>P\,,\qquad
m\geqslant m_0\,,
\end{equation}
in the case $p=n.$ If $p>n,$ then
\begin{equation}\label{eq1T}
M_p(\Gamma(|E^{m\,*}|, |\beta_m|, D\setminus \{x_0\}))\geqslant
P^{\frac{p}{n}}\cdot (m(D))^{\frac{n-p}{n}}\,,\qquad m\geqslant
m_0\,.
\end{equation}
Let us show that the relations~(\ref{eq11O}) and~(\ref{eq1T}) for
$p=n$ and $p>n,$ respectively, are impossible (in particular, each
of them contradicts definition of mapping~$f$
in~(\ref{eq2*A})--(\ref{eqA2})). Really, reasoning in the same way
as in case~1), by the ratio~(\ref{eq16A}) and~(\ref{eq17A}) we will
have that $f(\alpha_m)\subset B(z_1, r_1)$ and
$|D^{m\,}|=|f(D^{m\,*})|\subset {\Bbb R}^n\setminus \overline{B(z_1,
r_2)},$ where
$$z_1=f(\gamma^1_{k_0}(t_m))\,,$$
\begin{equation}\label{eq20C}
r_1=|a-z_1|+\varepsilon_0\,, \qquad r_2=|b-z_1|-\varepsilon_0\,.
\end{equation}
Then, similarly to~(\ref{eq5D}), we have:
\begin{equation}\label{eq5F}
\Gamma(E^{m\,*}, |\beta_m|, D\setminus\{x_0\})>\Gamma_f(z_1, r_1,
r_2)\,.
\end{equation}
In turn, by~(\ref{eq5F}) we have the following:
$$M_p(\Gamma(E^{m\,*}, |\beta_m|, D\setminus\{x_0\}))\leqslant$$
\begin{equation}\label{eq11R}
\leqslant M_p(\Gamma_f(z_1, r_1, r_2))\leqslant \int\limits_{A}
Q(y)\cdot \eta^p (|y-z_1|)\, dm(y)\,,
\end{equation}
where $A=A(z_1, r_1, r_2)$ and $\eta$ is an arbitrary Lebesgue
measurable function satisfying the relation~(\ref{eqA2}).
Now let us put
$$\eta(t)= \left\{
\begin{array}{rr}
\frac{1}{r_2-r_1}, & t\in [r_1, r_2],\\
0, & t\not\in [r_1, r_2]\,.
\end{array}
\right. $$
Note that, $\eta$ satisfies the relation~(\ref{eqA2}). Then
\begin{equation}\label{eq14*****A}
M_p(\Gamma(E^{m\,*}, |\beta_m|, D\setminus\{x_0\}))\leqslant
M_p(\Gamma_{f}(z^1, r_1, r_2))\leqslant
\frac{1}{(r_2-r_1)^p}\int\limits_{D^{\,\prime}}
Q(y)\,dm(y)<\infty\,.
\end{equation}
The relation~(\ref{eq14*****A}) contradicts with~(\ref{eq11O}) for
$p=n$ and~(\ref{eq1T}) for $p>n$. The resulting contradiction
completes the consideration case~3).

\medskip
{\bf Case~4)} is considered similarly to all the previous ones,
therefore its detailed consideration is left to the reader.~$\Box$
\end{proof}

\medskip
\section{Proof of the main results}

{\bf Proof of Theorem~\ref{th1}.} Due to the compactness of the
space $\overline{{\Bbb R}^n},$ the set $C(x_0, f)$ is not empty.
Moreover, by Proposition~\ref{pr1} $C(x_0, f)$ is a continuum in
$\overline{{\Bbb R}^n}.$ Assume the opposite: let the continuum
$C(x_0, f)$ is nondegenerate. Then, if $C(x_0, f)\cap
D^{\,\prime}\ne\varnothing,$ then by the connectedness of $C(x_0,
f)$ and the openness of $D^{\,\prime}$ the set $C(x_0, f)\cap
D^{\,\prime}$ must also consist of an infinite number of points. By
Lemma~\ref{lem4}, the latter is not possible. Then $C(x_0, f)\subset
\partial D^{\,\prime},$ but this contradicts with Lemma~\ref{lem1}
and Proposition~\ref{pr5}, because $C(x_0, f)$ is non-degenerate.
The obtained contradiction indicates that the continuum $C(x_0, f)$
is degenerate. This means the possibility of continuous extension of
the mapping $f$ to the point~$x_0.$

\medskip The openness and discreteness of the extended mapping
$\overline{f}$ in $D$ follows from Proposition~\ref{pr4}. If $p=n$
and $\overline{f}(x_0)\ne \infty,$ the estimate~(\ref{eq2C}) follows
from Proposition~\ref{pr6}.~$\Box$

\medskip
{\bf The proof of Theorem~\ref{th2}} is similar to the proof of
Theorem~\ref{th1} (let us to demonstrate it). Due to the compactness
of the space $\overline{{\Bbb R}^n},$ the set $C(x_0, f)$ is not
empty. Moreover, by Proposition~\ref{pr1} $C(x_0, f)$ is a continuum
in $\overline{{\Bbb R}^n}.$ Assume the opposite: let the continuum
$C(x_0, f)$ be non-degenerate. Then, if $C(x_0, f)\cap
D^{\,\prime}\ne\varnothing,$ by the connectedness of $C(x_0, f)$ and
the openness of $D^{\,\prime},$ the set $C(x_0, f)\cap D^{\,\prime}$
must also consist of an infinite number of points. It is not
possible by Lemma~\ref{lem3}. Then $C(x_0, f)\subset
\partial D^{\,\prime}$ but this contradicts Lemma~\ref{lem1} and Proposition~\ref{pr5},
because $C(x_0, f)$ is non-degenerate. The obtained contradiction
indicates that the continuum $C(x_0, f)$ is degenerate. The latter
means the possibility of continuous extension of the mapping of $f$
to the point~$x_0.$

\medskip The openness and discreteness of the extended mapping
$\overline{f}$ in $D$ follows by Proposition~\ref{pr4}.~$\Box$

\section{Some examples}

\begin{example}\label{ex1}
First, consider a plane mapping $f(z)=z^m,$ $z_0=0,$ for some fixed
$m=1,2, \ldots.$ Let $z_0=0$ and $D={\Bbb D}=\{z\in {\Bbb C}:
|z|<1\}.$ This mapping is quasiregular and it may be considered in
the form $f:{\Bbb D}\setminus\{0\}\rightarrow {\Bbb C}.$ In
addition, $f$ satisfies the condition~(\ref{eq2*A}) for $p=2$ at any
point $y_0\in \overline{{\Bbb D}}$ with $Q(y)=m$ (see Theorem~3.2
in~\cite{MRV$_1$}, Remark~\ref{rem2} and comments made in
introduction). The mapping $f$ satisfies all conditions of
Theorems~\ref{th1} and~\ref{th2}. By these theorems, it has a
continuous extension to the point $z_0=0.$ Note that, the mapping
$f$ satisfies the condition $C(0, f)\subset \partial f({\Bbb
D}\setminus 0)=\partial({\Bbb D}\setminus \{0\}).$ Therefore, its
continuous extension to the point $z_0=0$ may be obtained also from
Corollary~\ref{cor2} (see also \cite[Theorem~1]{Sev$_2$}
and~\cite[Theorem~6]{Sev$_1$}).
\end{example}

\medskip
\begin{example}\label{ex2}
Now, consider the same mapping, but at some point $z_0\ne 0.$ Let
$z_0=re^{i\varphi},$ $0\leqslant r<\infty,$ $\varphi\in[0, 2\pi).$
Consider $r_0<|z_0|$ and the domain $D:=B(z_0, r_0).$ Consider the
mapping $f:D\setminus\{z_0\}\rightarrow{\Bbb C}.$ Then
$z^m_0=r^me^{im\varphi}.$ Let $m$ be so large that the point
$z_1:=re^{i(\varphi+\frac{2\pi}{m})}$ still belongs to the domain
$D$ and is not equal to $z_0.$ Then $z_1^m=z_0^m=f(z_0)=f(z_1).$
Therefore, $C(z_0, f)\subset f({\Bbb D}\setminus 0).$ In this case,
the mapping $f$ is not closed, therefore the previous results
(proved earlier by the second co-author) cannot be applied to it
(see \cite[Theorem~1]{Sev$_2$} and~\cite[Theorem~6]{Sev$_1$}).
However, Theorems~\ref{th1} and~\ref{th2} are applicable; this
mapping is quasiregular, therefore $f$ satisfies the
condition~(\ref{eq2*A}) for $p=2$ at any point $y_0\in \overline{D}$
with $Q(y)=m$ (see Theorem~3.2 in~\cite{MRV$_1$}, Remark~\ref{rem2}
and comments, made in the introduction). The function $Q(y)=m$ is
integrable in $f(D\setminus \{z_0\}),$ because it is continuous, and
the closure of the set $f(D\setminus \{z_0\})$ is compact in ${\Bbb
C}.$ Of course, $f$ is open and discrete (see \cite[Theorems~6.3.II
and 6.4.II]{Re}).
\end{example}

\medskip
\begin{example}\label{ex3}
Let $n\geqslant 2,$ $m\in {\Bbb N},$ $x\in {\Bbb B}^n,$ and let
$f(x)=(r\cos m\varphi, r\sin m\varphi, x_3,\ldots, x_n),$ where
$x=(x_1,x_2,\ldots, x_n),$ $r=\sqrt{x_1^2+x_2^2}.$ For this mapping
$K_O(x, f)=m^{n-1}$ (see item~4 Ch.~I in~\cite{Re}). Obviously,
$N(f, {\Bbb B}^n)=m.$ Assuming that $x_0=0$ and consider mapping
$f:{\Bbb B}^n\setminus\{0\}\rightarrow {\Bbb R}^n,$ we may see that
$f$ satisfies the condition~(\ref{eq2*A}) for $p=n$ at any point
$y_0\in \overline{{\Bbb B}^n}$ with $Q(y)=m$ (see Theorem~3.2
in~\cite{MRV$_1$}, Remark~\ref{rem2} and comments made in the
introduction). The mapping $f$ satisfies all conditions of
Theorems~\ref{th1} and~\ref{th2}. According to these theorems, it
has a continuous extension to the point $x_0=0.$ Reasoning similarly
to the above Example~\ref{ex2} it is possible to construct the ball
$U=B(x_0, r_0)\subset {\Bbb B}^n,$ $x_0\ne 0,$ in which the
corresponding mapping $f:U\setminus\{x_0\}\rightarrow {\Bbb R}^n$ is
not a homeomorphism and satisfies all the conditions of
Theorems~\ref{th1} and~\ref{th2} for $Q(y)=m$.
\end{example}

\medskip
\begin{example}\label{ex4}
Let $x_0\in {\Bbb R}^n,$ $r_0>0,$ $n\geqslant 2.$ Put
$f(x)=\frac{r_0x}{|x|\log\frac{e}{|x|}}+x_0,$ $x\in B(x_0, 1),$
$f(0)=x_0.$ Let us put $g(y):=f^{\,-1}(y).$ Then $g$ is defined in
the ball $B(x_0, r_0)$ and $g(B(x_0, r_0))={\Bbb B}^n.$ Reasoning
similarly to \cite[Proposition~6.3]{MRSY}, it may be shown that $g$
satisfies the relations~(\ref{eq2*A})--(\ref{eqA2}) at any point
$y_0\in \overline{{\Bbb B}^n}$ for
$Q=Q(y)=\log^{n-1}\left(\frac{e}{|y|}\right).$

Note that $Q\in L^1({\Bbb B}^n).$ Indeed, by the Fubini theorem, we
assume that
$$\int\limits_{{\Bbb B}^n}Q(y)\,dm(y)=\int\limits_0^1\int\limits_{S(0, r)}
\log^{n-1}\left(\frac{e}{|y|}\right)d\mathcal{H}^{n-1}(y)dr=$$
$$=\omega_{n-1}\int\limits_0^1
r^{n-1}\log^{n-1}\left(\frac{e}{r}\right)\,dr\leqslant
\omega_{n-1}e^{n-1}\int\limits_0^1dr=\omega_{n-1}e^{n-1}<\infty\,,$$
where $\omega_{n-1}$ denotes the area of the unit sphere ${\Bbb
S}^{n-1}$ in ${\Bbb R}^n.$ Therefore, the mapping $g:B(x_0,
r_0)\setminus\{x_0\}\rightarrow {\Bbb R}^n$ satisfies all conditions
of Theorems~\ref{th1} and~\ref{th2}. The mapping $g$ is a
homeomorphism, and due to any of these theorems, it has a continuous
extension to a point~$x_0.$
\end{example}

\medskip
\begin{example}\label{ex5}
Let $f$ be one of the mappings from Examples~\ref{ex1}--\ref{ex3}.
Put $x_0=0$ and let $g_0:{\Bbb B}^n\setminus\{0\}\rightarrow {\Bbb
R}^n$ be a mapping from Example~\ref{ex4} for $x_0=0.$ Put
$$G_0(x):=(g_0\circ f)(x)\,.$$
In terms of the definition~(\ref{eq2*A}), it may be shown that
\begin{equation}\label{eq25}
f(\Gamma_{G_0}(y_0, r_1, r_2))\subset \Gamma_{g_0}(y_0, r_1, r_2)\,.
\end{equation}
Then, by~(\ref{eq25}) and Remark~\ref{rem2}, we will have that
$$M(\Gamma_{G_0}(y_0, r_1, r_2))\leqslant M_0\cdot f(\Gamma_{G_0}(y_0, r_1, r_2))
\leqslant$$
\begin{equation}\label{eq26}
\leqslant M_0\cdot M(\Gamma_{g_0}(y_0, r_1, r_2))\leqslant M_0\cdot
\int\limits_{({\Bbb B}^n\setminus\{0\})\cap A(y_0,r_1,r_2)}
Q(y)\cdot \eta^n (|y-y_0|)\, dm(y)\,,
\end{equation}
where $\eta: (r_1,r_2)\rightarrow [0,\infty ]$ is arbitrary Lebesgue
measurable function such that
$$
\int\limits_{r_1}^{r_2}\eta(r)\, dr\geqslant 1\,,
$$
and $Q=Q(y)=\log^{n-1}\left(\frac{e}{|y|}\right).$ (Here
$M_0\geqslant 1$ is equal to some constant, its own for each
Examples~\ref{ex1}--\ref{ex3}). From~(\ref{eq26}) it follows that,
the mapping $G_0:{\Bbb B}^n\setminus\{0\}\rightarrow {\Bbb R}^n$
satisfies all conditions of Theorems~\ref{th1} and~\ref{th2}. This
mapping is open and discrete as a superposition of a quasiregular
mapping $f$ with the homeomorphism $g_0,$ so that by any of
Theorems~\ref{th1} or~\ref{th2} it has a continuous extension to the
origin. Note that, $C(0, G_0)\subset \partial ({\Bbb
B}^n\setminus\{0\}).$ \end{example}

\medskip
\begin{example}\label{ex6}
Finally, let $f$ be one of the mappings from
Examples~\ref{ex1}--\ref{ex3}. Let $x_0\ne 0$ and let $B(x_0, r_0)$
be some ball in which mapping $f$ is not a homeomorphism (this may
be achieved by increasing of the index $m$). Let $f(x_0)=p_0.$
Choose $R_0>0$ is so large that $f(B(x_0, r_0))\subset B(p_0, R_0).$

\medskip
Put $F(x)=\frac{R_0x}{|x|\log\frac{e}{|x|}}+p_0,$ $x\in B(x_0, 1),$
$f(0)=p_0.$ Let $G(y):=F^{\,-1}(y)$ and
$$G_1(x):=(G\circ f)(x)\,.$$
Reasoning similarly to Example~\ref{ex5}, it may be shown that the
mapping $G_1$ satisfies the conditions~(\ref{eq2*A})--(\ref{eqA2})
for $p=n$ with $Q=Q(y)=\log^{n-1}\left(\frac{e}{|y|}\right).$ (Here
$M_0\geqslant 1$ is equal to some constant, its own for each of the
Examples~\ref{ex1}--\ref{ex3}). From here, it follows that, the
mapping $G_1:B(x_0, r_0)\setminus\{x_0\}\rightarrow {\Bbb R}^n$
satisfies all conditions of Theorems~\ref{th1} and~\ref{th2}. This
mapping is open and discrete as a superposition of some quasiregular
mapping $f$ with some homeomorphism $G,$ and by any of
Theorems~\ref{th1} or~\ref{th2} it has a continuous extension to the
point~$x_0.$ Note that $C(0, G_1)\subset (G_1(B(x_0,
r_0)\setminus\{x_0\})).$ \end{example}

\medskip
{\bf \noindent Victoria Desyatka} \\
Zhytomyr Ivan Franko State University,  \\
40 Bol'shaya Berdichevskaya Str., 10 008  Zhytomyr, UKRAINE \\
victoriazehrer@gmail.com

\medskip
\medskip
{\bf \noindent Evgeny Sevost'yanov} \\
{\bf 1.} Zhytomyr Ivan Franko State University,  \\
40 Bol'shaya Berdichevskaya Str., 10 008  Zhytomyr, UKRAINE \\
{\bf 2.} Institute of Applied Mathematics and Mechanics\\
of NAS of Ukraine, \\
19 Henerala Batyuka Str., 84 116 Slavyansk,  UKRAINE\\
esevostyanov2009@gmail.com

\end{document}